\documentclass[10pt]{article}
\usepackage[margin=1.3in]{geometry}

\usepackage{hyperref}
\newcommand{\email}[1]{\href{mailto:#1}{#1}}

\usepackage{amsmath,amssymb,amsfonts}
\usepackage{dirtytalk} 

\usepackage{graphicx}
\usepackage{subfig} 
\usepackage{tikz}
\usetikzlibrary{patterns,intersections,arrows.meta,decorations.markings,cd,spy}
\usepackage{tikz-3dplot}
\usepackage{epstopdf}

\usepackage[noabbrev]{cleveref}
\crefname{equation}{}{}
\numberwithin{equation}{section}

\usepackage{amsthm}

\newtheorem{lemma}{Lemma}[section]
\newtheorem{proposition}{Proposition}[section]

\theoremstyle{definition}
\newtheorem{definition}{Definition}[section]
\newtheorem{remark}{Remark}[section]
\newtheorem{example}{Example}[section]

\DeclareMathOperator{\curl}{curl}
\DeclareMathOperator{\sign}{sign}
\DeclareMathOperator{\diag}{diag}

\DeclareMathOperator*{\argmin}{argmin}

\newcommand{\calG}{\mathcal{G}}

\newcommand{\calP}{\mathcal{P}}
\newcommand{\calB}{\mathcal{B}}
\newcommand{\calV}{\mathcal{V}}
\newcommand{\calM}{\mathcal{M}}
\newcommand{\calX}{\mathcal{X}}
\newcommand{\calY}{\mathcal{Y}}

\newcommand{\x}{\boldsymbol{x}}
\newcommand{\y}{\boldsymbol{y}}
\newcommand{\m}{\boldsymbol{m}}
\newcommand{\p}{\boldsymbol{p}}
\newcommand{\pphi}{\boldsymbol{\varphi}}

\newcommand{\PP}{\boldsymbol{P}}
\newcommand{\bcdot}{\boldsymbol{\cdot}}
\newcommand{\bmu}{\boldsymbol{\mu}}

\newcommand{\D}{\mathrm{D}}
\newcommand{\wrom}{\mathrm{w}}
\newcommand{\nrom}{\mathrm{n}}
\newcommand{\erom}{\mathrm{e}}
\newcommand{\srom}{\mathrm{s}}
\newcommand{\tL}{\mathrm{L}}
\newcommand{\tR}{\mathrm{R}}
\newcommand{\I}{\mathrm{int}}

\newcommand{\dx}{\,\mathrm{d}\x}
\newcommand{\dy}{\,\mathrm{d}\y}

\newcommand{\ds}{\,\mathrm{d}s}

\newcommand{\unitv}[1]{\hat{\boldsymbol{#1}}}

\newcommand{\on}{\text{on } \; }
\newcommand{\RI}{\mathbb{R}}
\newcommand{\GamintX}{\Gamma_\I^\calX}
\newcommand{\diff}[2]{\frac{\partial #1}{ \partial #2}}

\newcommand{\overbar}[1]{\mkern 1.5mu\overline{\mkern-4mu#1\mkern-1.5mu}\mkern 1.5mu}

\newcommand{\expseq}[2]{^{#1,(#2)}}

\newcommand{\littletaller}{\mathchoice{\vphantom{\big|}}{}{}{}}
\newcommand\restr[2]{{
  \left.\kern-\nulldelimiterspace 
  #1 
  \littletaller 
  \right|_{#2} 
  }}

\title{A Mixed Elliptic--Hyperbolic Monge--Ampère Problem: Modeling and a Least-Squares Solver}
\author{G.V. Carvalho\thanks{CASA, Department of Mathematics and Computer Science, Eindhoven University of
Technology, PO Box 513, 5600 MB Eindhoven, The Netherlands 
(\email{g.vidal.de.carvalho@tue.nl},
\email{m.j.h.anthonissen@tue.nl}, 
\email{j.h.m.tenthijeboonkkamp@tue.nl}, 
\email{l.kusch@tue.nl},
\email{k.mitra@tue.nl}.
\email{w.l.ijzerman@tue.nl}), }
, M.J.H. Anthonissen\footnotemark[1]
, J.H.M ten Thije Boonkkamp\footnotemark[1],
\\
L.~Kusch\footnotemark[1]
, K.~Mitra\footnotemark[1]
, W.L. IJzerman\footnotemark[1] \thanks{Signify Research, High Tech Campus 7, 5656 AE Eindhoven, The Netherlands (\email{wilbert.ijzerman@signify.com})}
}

\begin{document}

\maketitle

\begin{abstract}
We consider a transport problem motivated by freeform optical design of a parallel-to-far-field reflector. Classical formulations rely on either elliptic or hyperbolic Monge--Ampère equations, leading to globally convex/concave or saddle-shape reflector geometries, respectively. However, the design of smooth reflector arrays may require optical surfaces combining both types of curvature. To address this, we introduce a mixed elliptic--hyperbolic Monge--Ampère system obtained by partitioning the source domain into one elliptic and one hyperbolic subregion, and we develop a least-squares method for its numerical approximation. Numerical examples demonstrate the ability of the proposed framework to compute both isolated and periodic mixed-type transport solutions, yielding smooth reflector configurations.
\end{abstract}

\textbf{Keywords.} Mixed elliptic-hyperbolic equation, iterative least-squares algorithm, freeform optical design.

\textbf{MSC codes.} 35J96, 35M10

\section{Introduction}

The Monge-Ampère (MA) equation is a non-linear partial differential equation arising in various contexts such as prescribed curvature problems, optimal transport and freeform optical design. The present work is motivated by the latter, namely the design of freeform parallel-to-far-field reflectors.

In parallel-to-far-field reflector design we consider a beam of parallel light rays being emitted from a source region $\calX \subset \RI^2$ (typically a bounded open domain) with light density $f$, and a target region $\mathcal T \subset S^2$ (representing a set of reflected directions) with light density $g$, where $S^2 \subset \RI^3$ is the unit 2-sphere. The objective is to obtain a reflector surface parametrized by a function $u:\calX \to \RI$ (representing the height of the reflector at each point), which redirects the emitted light onto $\mathcal T$ while inducing the desired light intensity $g$. It is useful to parametrize $\mathcal{T}$ in stereographic coordinates and consequently we then refer to the target domain as $\calY \subset \RI^2$. 
Mathematically, this can be formulated as finding a conservative mapping $\m:\calX \to \calY$ satisfying
\begin{equation}\label{eq:transport_problem}
    \int_B g\dy = \int_{\m^{-1}(B)} f\dx, \quad \m = \nabla u,
\end{equation}
for every $B \subset \calY$. Throughout the paper the function $u$ will be denoted as the \textit{potential function}.
The first equation is the usual push-forward of measures, and the second equation guarantees that a reflector function can be recovered from $\m$.
For more details see for example \cite{jan_frontiers} for a recent survey on inverse freeform design in non-imaging optics. 

This problem has been mainly addressed via optimal transport techniques, leading to the second boundary value problem for the elliptic Monge--Ampère equation:
\begin{subequations}
\label{sys:MA_std}
\begin{align}
    \label{eq:MA_std}
    \det \D^2 u &= \eta \frac{f}{g \circ \nabla u} \quad \on \calX,
    \\
    \label{eq:trasnport_BC}
    \nabla u(\calX) &= \calY,
\end{align}
\end{subequations}
for $\eta=1$. This approach is advantageous since Brenier \cite{Brenier} guarantees existence and uniqueness of convex solutions to \cref{sys:MA_std}. There are also interior and boundary regularity results, namely Caffarelli \cite{caf_interior_reg,caf_boundary1,caf_boundary2} and Urbas \cite{MR1454261}. Numerical schemes for \cref{sys:MA_std} include Froese \cite{froese} and a least-squares method initially proposed by Caboussat et al. \cite{caboussat} for the Dirichlet problem, and later adapted by Prins et al. \cite{corien_paper} to \cref{sys:MA_std}.

For $\eta = -1$, equation \cref{eq:MA_std} becomes a hyperbolic MA equation. Compared to its elliptic counterpart, the hyperbolic system \cref{sys:MA_std} remains poorly understood, with fundamental questions concerning existence and uniqueness still open. Regarding numerical methods, recently Bertens et al. \cite{maikel} adapted the least-squares method developed in \cite{corien_paper} to the hyperbolic system \cref{sys:MA_std}.

\begin{figure}[ht]
    \centering
    \hspace{-4mm}
    \subfloat[Periodic array of convex elements exhibiting non-smooth interfaces. Note the sharp cusps between elements.]
    {\label{fig:array_convex}
        \centering
        \includegraphics[width=0.4\textwidth]{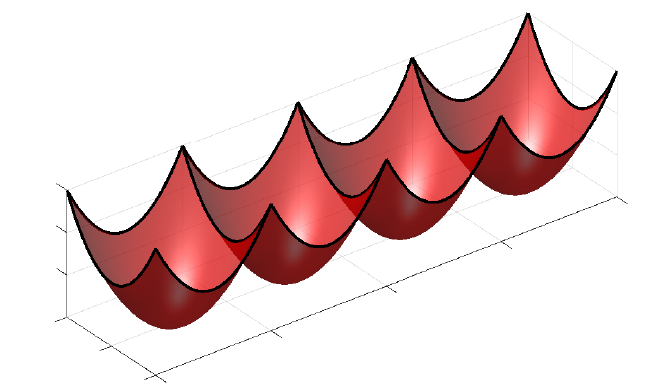}
    }
    \hspace{4mm}
    \subfloat[Smooth periodic array combining convex and saddle-like elements.]
    {\label{fig:array_smooth}
        \centering
        \includegraphics[width=0.4\textwidth]{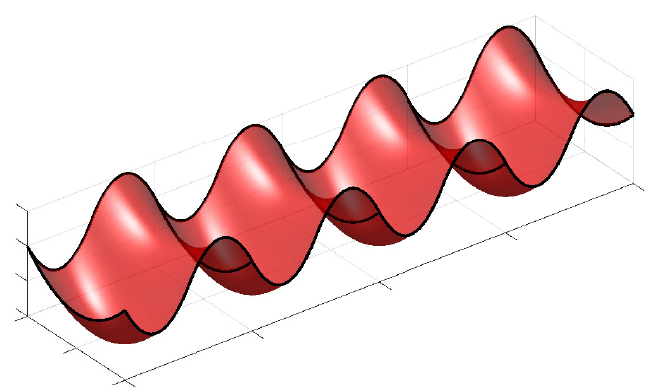}
    }
    \caption{Comparison between non-smooth and smooth periodic arrays.}
    \label{fig:periodic_arrays}
\end{figure}

In practical illumination applications, a single optical element is often replicated to form reflector or lens arrays, thereby increasing the total output intensity while preserving the desired target light illumination profile. When such arrays are constructed exclusively from convex or concave elements, the resulting surfaces exhibit discontinuities in slope along the interfaces between neighboring elements (see \cref{fig:array_convex}). These cusps, however, may generate undesirable optical artifacts
if they aren not manufactured perfectly sharp.
In this work we look into the design of smooth periodic arrays, which requires the consideration of surfaces with regions of different curvature signatures, by combining locally convex/concave functions with saddle-shape functions (see \cref{fig:array_smooth}).

Motivated by this application we introduce a novel mixed Monge--Ampère formulation of the transport problem \cref{eq:transport_problem}. The source domain is split into two subregions, a \textit{left} (denoted as $\calX^\tL$) and a \textit{right} (denoted as $\calX^\tR$) subregion, on which the potential function $u$ exhibits different curvature signatures. It is convex on one portion and saddle-shape on the other which leads to a coupled system consisting of an elliptic and a hyperbolic MA equation. 
In this formulation, each point in the target $\calY$ is illuminated by two points from the source, one from each subregion. Consequently the full target density has two separate density contributions, $g^\tL$ and $g^\tR$, which correspond to the reflected light from each subregion.
Ultimately, to satisfy the transport problem \cref{eq:transport_problem} we enforce the coupling equation $g^\tL+g^\tR = g$.
We found that the resulting system is underdetermined and may admit multiple solutions, but rather than addressing analytical questions of existence and regularity, the focus of this work is on the formulation of the mixed system and the development of a numerical method.

Mixed-type equations arise in several contexts, most notably in connection with transonic flow problems, where the Tricomi equation is the classical example. Particularly for the Monge--Ampère equation, mixed elliptic--hyperbolic equations have been studied for isometric embedding problems (see for example \cite{mixMA}) but to the best of the authors' knowledge, mixed Monge--Ampère formulations associated with mass transport problems have not been previously considered.

The contributions of this paper are threefold. First, we propose a novel mixed elliptic–hyperbolic Monge–Ampère formulation for a transport problem arising in optical design. Second, we develop a numerical scheme for approximating solutions of this mixed system. Third, we demonstrate the ability of the proposed approach to generate smooth reflector configurations, including periodic solutions.

This article is organized as follows. \Cref{sec:properties_map} summarizes relevant properties of solutions to the Monge--Ampère equation \cref{sys:MA_std}. The mixed elliptic–hyperbolic Monge–Ampère formulation is introduced in \Cref{sec:Mixed_MA}. In \Cref{sec:mix_LSA} we introduce a least-squares method for the mixed system, and ultimately, the numerical results are presented in \Cref{sec:results}.

\section{Properties of mapping solutions}\label{sec:properties_map}

In this section we briefly discuss some properties of $\m = \nabla u$ for solutions of \cref{sys:MA_std}. We will be assuming in this paper that the sets $\calX$ and $\calY$ are axis-aligned rectangular sets and without loss of generality, by a translation argument, we can assume them to be centered around $0$, that is,
\begin{equation}\label{def:XY}
    \calX = (-a, a) \times (-b,b), \quad \calY = (-c,c) \times (-d,d)
\end{equation}
for some $a,b,c,d>0$. 
We note however, that some results in this section hold for more general domains.

Furthermore, we will consider mass densities $f\in L^2(\calX)$ and $g \in L^2(\calY)$ bounded away from zero and infinity,
\begin{equation}\label{eq:bounded_away}
    c<f(\x),g(\y) < C, \quad \x\in\calX,\;\y\in\calY,
\end{equation}
for $c,C\in(0,+\infty)$ with $c<C$. These densities must have equal mass and we suppose without loss of generality that they have mass equal to $1$, i.e.,
\begin{equation}\label{eq:total_energy_conservation}
    \int_\calX f\dx = \int _\calY g\dy = 1.
\end{equation}

We proceed to write system \cref{sys:MA_std} in terms of the \emph{mapping} function $\m = (m_1,m_2) := \nabla u$. Here we are assuming $u$ to be a classical solution of \cref{sys:MA_std}, that is $u\in C^2(\calX)$. The function $u$ will often be referred as the \emph{potential} function of $\m$. Consequently, from \cref{sys:MA_std} we obtain,
\begin{subequations}\label{sys:MA_map}
\begin{alignat}{2}\label{eq:MA_map}
     &\det \D\m = \eta \displaystyle \frac{f}{g\circ\m} \quad &&\on \calX,
 \\
 \label{eq:curl_map}
 \displaystyle & \curl \m = \mathbf 0, &&\on\calX,
 \\
 \label{eq:transport_boundary_map}
 &\m(\calX) = \calY,
\end{alignat}
\end{subequations}
for $\eta\in\{-1,1\}$. As with \cref{sys:MA_std}, we refer to this system as \textit{elliptic} if $\eta=1$ and as \textit{hyperbolic} if $\eta = -1$.
In \cref{eq:curl_map} we write  $\curl \m = \partial_{x_1}m_2 - \partial_{x_2}m_1$ as the two-dimensional curl.

As we see in the following proposition, solutions to \cref{sys:MA_map} are invertible. For the elliptic system this is immediately clear since $\m$ is the gradient of a strictly convex or concave function. For the hyperbolic system the proof is slightly more involved.
\begin{proposition}[Diffeomorphism]\label{prop:diffeo}
    Suppose $\m\in [C^1(\calX)]^2$ satisfies \cref{eq:MA_map,eq:transport_boundary_map} with $f,g$ satisfying \cref{eq:bounded_away,eq:total_energy_conservation}. Then $\m:\calX \to \calY$ is a diffeomorphism.
\end{proposition}
\begin{proof}
    Surjectivity follows directly from \cref{eq:transport_boundary_map}. To prove injectivity, we first compute
    \[
    \int_\calX f\dx = \int_\calX  (g\circ \m)\, |\det \D\m|\dx = \int_\calY \#\Big(\m^{-1}\big(\{\y\}\big) \Big) \, g(\y) \dy,
    \]
    where in the first equality we used \cref{eq:MA_map} and in the second we used the area formula (see \cite[Theorem~3.2.3]{federer1969geometric}), where  $\#\big(\m^{-1}\big(\{\y\}\big) \big)$ denotes the cardinality of the set $\m^{-1}\big(\{\y\}\big) $. From \cref{eq:total_energy_conservation} we conclude that $\m^{-1}\big(\{\y\}\big) = 1$ for $\y$ almost everywhere.
    
    Suppose now that there exist distinct $\x_0,\x_1 \in \calX$ such that $\m(\x_0)=\m(\x_1)$. By \cref{eq:MA_map}, $\m$ is a local diffeomorphism therefore we can find disjoint open sets $A_0,A_1 \subset \calX$ such that $\x_0 \in A_0$, $\x_1 \in A_1$ and $\m(A_0) = \m(A_1)$ which has non-negative measure since $\m$ is an open map, contradicting $\m^{-1}\big(\{\y\}\big) = 1$ for $\y$ a.e.. Therefore, $\m$ is a bijective local diffeomorphism, and therefore a global diffeomorphism.
\end{proof}

We look now at the transport condition \eqref{eq:transport_boundary_map}. This is a global constraint on $\calX$ but the following proposition shows that, for mappings satisfying \eqref{eq:MA_map}, it is sufficient to prescribe the behavior of the mapping on the boundary, namely $\m(\partial \calX) = \partial \calY$. This characterization will be particularly useful in the numerical treatment developed later.

\begin{proposition}[Transport boundary condition]
    \label{prop:transport_boundary}
    If $\m\in [C^1(\calX) \cap C^0(\overbar \calX)]^2$ satisfies \cref{eq:MA_map} with $f,g$ satisfying \cref{eq:bounded_away,eq:total_energy_conservation}, then 
    \begin{equation}
        \m(\calX) = \calY \quad \iff \quad \m(\partial\calX) = \partial \calY.
    \end{equation}
\end{proposition}
\begin{proof}
$(\implies)$
Since $\m$ is continuous on $\overline{\calX}$, $\overline{\calX}$ is compact, and
$\m(\calX)=\calY$, then we have
$\m(\overbar {\calX})=\overbar{\calY}$.
By Proposition~\ref{prop:diffeo}, $\m$ is open and injective on $\calX$. One can then check that $\m(\partial\calX)\cap \m(\calX)=\emptyset$, which leads to $\m(\partial \calX) = \partial \calY$.

$(\impliedby)$ By \cref{eq:MA_map}, $\m$ is an open map, hence $\m(\calX)$ is open, so $\m(\calX) \cap \partial \m(\calX) = \emptyset$.

We first show $\m(\calX)\subset \calY$. By contradiction suppose $\m(\calX) \setminus \overbar \calY \neq \emptyset $. Continuity of $\m$ and boundedness of $\calX$ imply the compactness of $\overbar{ \m(\calX)}$. Hence, let $\y_0\in\partial  \m(\calX) \setminus \overbar \calY$ (take for example the point in $\overbar{\m(\calX)}$ maximizing the distance to $\overbar \calY$). By continuity of $\m$ there exists $\x_0 \in \overbar\calX$ such that $\m(\x_0) = \y_0$. Since $\m(\calX) \cap \partial \m(\calX) = \emptyset$ and $\y_0\in \partial \m(\calX)$, then $\x_0 \in \partial \calX$, contradicting $\m(\partial \calX) = \partial\calY$.

Thus $\m(\calX)\subset \calY$. Note then that $\overbar{\m(\calX)} \cap \calY = \big( \m(\calX) \cup \m(\partial \calX) \big) \cap \calY = \m(\calX)$ so $\m(\calX) $ is open and closed in $\calY$. Since $\calY$ is connected and $\m(\calX) \neq \emptyset$ we get that $\m(\calX) = \calY$.
\end{proof}

\begin{remark}[Boundary regularity assumption]\label{rem:regularity_assumption}
    In the remainder of this paper we assume that solutions to \cref{sys:MA_map} satisfy $\m\in [C^1(\overbar \calX)]^2$. This regularity assumption together with \cref{eq:bounded_away,prop:diffeo} allows us to treat $\m:\overbar \calX \to \overbar \calY$ as a diffeomorphism between the closed sets. Consequently, this simplifies the presentation of the boundary analysis below. Although several of the arguments can be extended under weaker assumptions, requiring only continuity up to the boundary, this would be at the expense of additional technicalities that are beyond the scope of the present work.
\end{remark}

We now go through some properties of the boundary map $\restr{\m}{\partial \calX}:\partial \calX \to \partial \calY$. 
We define a \textit{simple closed curve} to be a continuous, non-intersecting curve $\boldsymbol{\gamma} : [0,1] \to \RI^2$ with $\boldsymbol{\gamma}(0) = \boldsymbol{\gamma}(1)$. By the Jordan Curve Theorem, $\boldsymbol{\gamma}([0,1])$ separates $\RI^2$ in two disjoint connected components, an interior and an exterior component. We say that a simple closed curve has a positive (negative) orientation if the interior of the curve is always to the left (right) when traveling through $\boldsymbol{\gamma}$. In particular, if $\boldsymbol\gamma $ is regular, meaning that $\boldsymbol\gamma' \neq \boldsymbol{0}$, then orientation of $\boldsymbol\gamma$ is equivalent to the orientation of the basis $(\unitv{n}, \boldsymbol{\gamma}')$ where $\unitv{n}$ is the outward normal to $\boldsymbol{\gamma}$. That is, $\boldsymbol{\gamma}$ is positively oriented if $\det (\unitv{n}, \boldsymbol{\gamma}') > 0$, and negatively oriented if $\det (\unitv{n}, \boldsymbol{\gamma}') < 0$.

\begin{lemma}[Orientation under a diffeomorphism]
    Let $\boldsymbol{\gamma} \in C^1([0,1])$ be a simple closed curve in $\calX$ with positive orientation satisfying $\boldsymbol{\gamma}' \neq \boldsymbol{0}$, and let $\m:\calX \to \calY$ be a diffeomorphism. Then $\m \circ \boldsymbol{\gamma}$ has positive orientation if $\det \D\m >0$ and negative orientation if $\det \D\m < 0$.
\end{lemma}
\begin{proof}
    Define $\unitv{n}:[0,1] \to S^1$ as the outward unit normal to $\boldsymbol{\gamma}$, and $\hat{\boldsymbol{\nu}}$ as the outward unit normal to $\m \circ \boldsymbol{\gamma}$. Since $\boldsymbol{\gamma}$ is positively oriented then $\det(\unitv{n},\boldsymbol{\gamma}') > 0$. Furthermore, since $\m$ is a diffeomorphism then the interior of the Jordan curve $\boldsymbol{\gamma}$ is also mapped to the interior of $\m \circ \boldsymbol{\gamma}$. In particular $(\D\m) \,  \unitv{n}$ points outwards of $\m \circ \boldsymbol{\gamma}$, therefore $((\D\m) \, \unitv{n} )\cdot \hat{\boldsymbol{\nu}} > 0$. Note that $(\m \circ \boldsymbol{\gamma})'=(\D\m) \, \boldsymbol{\gamma}'$, thus the orientation of $\m \circ \boldsymbol{\gamma}$ is given by the orientation of the basis $(\hat{\boldsymbol{\nu}},(\D\m) \, \boldsymbol{\gamma}')$, which can be computed as
    \[
    \sign \left( \det(\hat{\boldsymbol{\nu}},\D\m \, \boldsymbol{\gamma}') \right) = \sign \left( \det(\D\m \, \unitv{n},\D\m \, \boldsymbol{\gamma}') \right) = \sign (\det \D\m) \sign \left( \det(\unitv{n},\boldsymbol{\gamma}') \right),
    \]
    where in the first equality we used $\D\m\, \unitv{n} \cdot \hat{\boldsymbol{\nu}} > 0$.  Since $\det(\unitv{n},\boldsymbol{\gamma}') > 0$ we conclude that $\m \circ \boldsymbol{\gamma}$ has positive orientation if $\det \D\m >0$ and negative orientation if $\det \D\m < 0$.
\end{proof}

\begin{remark}[Boundary orientation]\label{rem:boundary_orientation}
    The previous lemma allows us to characterize the orientation of the boundary map $\restr{\m}{\partial \calX} : \partial\calX \to \partial \calY$. More concretely, let $\boldsymbol\gamma:[0,1] \to \partial \calX$ be a simple closed curve with positive orientation and $\boldsymbol\gamma([0,1]) = \partial \calX$, and $\m:\overbar\calX \to \overbar\calY$ be a diffeomorphism. Then $\m \circ \boldsymbol\gamma :[0,1] \to \partial \calY$ is a simple closed curve with positive orientation if $\det \D\m > 0$, and with negative orientation if $\det \D\m <0$. 
\end{remark}

By the regularity assumption in Remark \ref{rem:regularity_assumption} we obtain some additional regularity of the boundary map. 

\begin{lemma}[Correspondence of corner points]\label{lem:corner_points}
    Let $\m:\overbar \calX \to \overbar \calY$ be a diffeomorphism. Then $\m$ maps corner points of $\calX$ to corner points of $\calY$.
\end{lemma}
\begin{proof}
    From \cref{eq:bounded_away,eq:MA_map}, we can conclude that $|\det \D\m| \neq 0$ on $\partial \calX$. In particular, the image of a smooth segment of $\partial \calX$ must be a smooth curve. Therefore each corner point of $\partial \calY$ must be the image of a corner point of $\partial \calX$. Since $\m$ is injective, this corner map is one-to-one.
\end{proof}

\begin{figure}[ht]
    \centering
    \subfloat[Domain $\calX$.]{
        \label{subfig:Xseg}
        \centering 
        \includegraphics{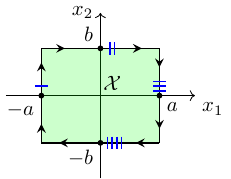}
    }
    \subfloat[Segmentation \cref{eq:bseg_1}]{
        \centering
        \label{subfig:bseg_1}
        \includegraphics{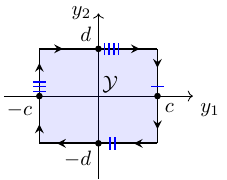}
    
    }
    \subfloat[Segmentation \cref{eq:bseg_2}]{
        \centering
        \label{subfig:bseg_2}
        \includegraphics{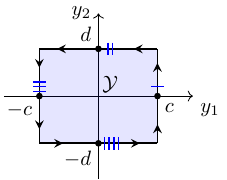}
    }

    \caption{Depiction of the boundary segmentations \cref{eq:bseg_1,eq:bseg_2}. In \cref{subfig:Xseg} we mark each segment of $\partial\calX$ with blue ticks. In \cref{subfig:bseg_1,subfig:bseg_2}, the marks of each segment of $\partial \calY$ follow from the boundary segmentations \cref{eq:bseg_1,eq:bseg_2}, respectively. The arrows emphasize the boundary orientation.}
    \label{fig:segmentation_single}
\end{figure}

This correspondence of corner points induces also a correspondence of the smooth boundary segments of $\partial \calX$ and $\partial\calY$.
Let $\Gamma_\wrom^\calX$, $\Gamma_\nrom^\calX$, $\Gamma_\erom^\calX$ and $\Gamma_\srom^\calX$ be the west (left), north (top), east (right) and south (bottom) segments of $\partial \calX$, that is,
\[
\Gamma^\calX_\wrom = \{-a\}\times (-b,b), \;  \Gamma^\calX_\nrom = (-a,a) \times \{b\}, \;
\Gamma^\calX_\erom = \{a\}\times (-b,b), \;  \Gamma^\calX_\srom = (-a,a) \times \{-b\}.
\]
We define similarly $\Gamma_\wrom^\calY$, $\Gamma_\nrom^\calY$, $\Gamma_\erom^\calY$ and $\Gamma_\srom^\calY$ as the boundary segments of $\partial \calY$. 
\begin{definition}[Boundary segmentation]\label{def:segmentation}
    Let 
    \[
    \boldsymbol{\Gamma}^\calX = \{\Gamma_\wrom^\calX,\Gamma_\nrom^\calX,
    \Gamma_\erom^\calX,\Gamma_\srom^\calX\},
    \quad
    \boldsymbol{\Gamma}^\calY = \{\Gamma_\wrom^\calY,\Gamma_\nrom^\calY,
    \Gamma_\erom^\calY,\Gamma_\srom^\calY\}.
    \]
    Then a \textbf{boundary segmentation} of a mapping solution $\m:\overbar\calX \to \overbar \calY$ to \cref{sys:MA_map} is the assignment
    \[
    \Gamma^\calX_\ast \in \boldsymbol{\Gamma}^\calX \;\; \mapsto \;\; \m(\Gamma^\calX_\ast) \in \boldsymbol{\Gamma}^\calY,
    \]
    with $\ast\in\{\wrom,\nrom,\erom,\srom \}$.
\end{definition}
From Remark \ref{rem:boundary_orientation} and Lemma \ref{lem:corner_points}, the boundary segment assignment is cyclical and completely determined by the image of one single segment together with the sign of $\det\D\m$. As an example, if $\m$ solves the elliptic system \cref{sys:MA_map},
\begin{equation}
\label{eq:bseg_1}
\m(\Gamma^\calX_\wrom) = \Gamma^\calY_\erom  \quad \implies \quad \m(\Gamma^\calX_\nrom) = \Gamma^\calY_\srom, \;\; 
\m(\Gamma^\calX_\erom) = \Gamma^\calY_\wrom, \;\; \m(\Gamma^\calX_\srom) = \Gamma^\calY_\nrom.
\end{equation}
Note that the image of traveling clockwise along $\partial\calX$ remains clockwise along $\partial\calY$ (see \cref{subfig:Xseg,subfig:bseg_1}).

On the contrary, if $\m$ is a solution to the hyperbolic system \cref{sys:MA_map}, then,
\begin{equation}
\label{eq:bseg_2}
\m(\Gamma^\calX_\wrom) = \Gamma^\calY_\erom  \quad \implies \quad 
\m(\Gamma^\calX_\nrom) = \Gamma^\calY_\nrom, \; 
\m(\Gamma^\calX_\erom) = \Gamma^\calY_\wrom, \; \m(\Gamma^\calX_\srom) = \Gamma^\calY_\srom.
\end{equation}
Note now that the image of traveling clockwise along $\partial\calX$ becomes counter-clockwise along $\partial\calY$, as depicted in \cref{subfig:Xseg,subfig:bseg_2}.

We see that a priori, by Remark \ref{rem:boundary_orientation}, there are four possible boundary segmentations for solutions to the elliptic system, and four other for the hyperbolic system. But in fact, as stated in the following proposition, we can deem inadmissible two segmentations for the elliptic system.

\begin{proposition}[Elliptic segmentations]\label{prop:ell_seg}
    Let $\m\in [C^1(\overbar \calX)]^2$ be a solution to elliptic system \cref{sys:MA_map} (i.e. for $\eta = 1$). Then the boundary segmentation of $\m$ is 
    either\begin{equation}
        \label{subeq:segmentation_convex}
        \m(\Gamma^\calX_\wrom) = \Gamma^\calY_\wrom, \quad
        \m(\Gamma^\calX_\nrom) = \Gamma^\calY_\nrom, \quad
        \m(\Gamma^\calX_\erom) = \Gamma^\calY_\erom, \quad
        \m(\Gamma^\calX_\srom) = \Gamma^\calY_\srom,
    \end{equation}
    implying that $\m=\nabla u$ for a \textbf{convex} $u$, or,
    \begin{equation}
        \label{subeq:segmentation_concave}
        \m(\Gamma^\calX_\wrom) = \Gamma^\calY_\erom, \quad
        \m(\Gamma^\calX_\nrom) = \Gamma^\calY_\srom, \quad
        \m(\Gamma^\calX_\erom) = \Gamma^\calY_\wrom, \quad
        \m(\Gamma^\calX_\srom) = \Gamma^\calY_\nrom,
    \end{equation}
    implying that $\m = \nabla u$ for a \textbf{concave} $u$.
\end{proposition}
\begin{proof}
    If a horizontal segment is mapped to a vertical segment, then $m_1$ is constant along that segment, hence $\partial_{x_1} m_1 =0$. From \cref{eq:curl_map}, $\det \D\m = -(\partial_{x_1}m_2)^2\leq 0$ on this segment, which contradicts \cref{eq:MA_map} via \cref{eq:bounded_away}. A similar reasoning also implies that no vertical segment of $\partial\calX$ can be mapped to a horizontal segment of $\partial \calY$. 

    This leaves us only with \cref{subeq:segmentation_convex,subeq:segmentation_concave} as the possible boundary segmentations for the elliptic system \cref{sys:MA_map} as a consequence of positive orientation (Remark \ref{rem:boundary_orientation}). In the case of \cref{subeq:segmentation_convex} we see that 
    \[
    \partial_{x_2}m_1 = 0,\quad  \partial_{x_2}m_2>0, \quad \on \Gamma_\wrom^\calX.
    \]
    From \cref{eq:MA_map} and \cref{eq:curl_map} we conclude that $\D \m$ is symmetric positive definite (SPD) on $\Gamma^\calX_\wrom$. Because the Jacobian determinant of $\m$ never vanishes and $\D \m$ is always symmetric (cf. \cref{eq:curl_map}), its eigenvalues can vary continuously but can never cross zero. Since $\calX$ is connected, $\D\m$ remains SPD in $\calX$, hence $\m = \nabla u$ for convex $u$.
    
    An analogous argument can be applied to the case \cref{subeq:segmentation_concave}.
\end{proof}

\section{Mixed elliptic-hyperbolic Monge--Ampère system}\label{sec:Mixed_MA}

In this section we introduce a formulation of the transport problem \cref{eq:transport_problem} on rectangular domains that only admits solutions $\m = \nabla u$ where $u$ is convex on a left subdomain and saddle-like on the remaining subdomain, while preserving differentiability of $u$ in the entire domain. The graph of $u$ will then resemble an element of the surface in  \cref{fig:array_smooth}.

\begin{figure}[ht]
    \centering
    \includegraphics{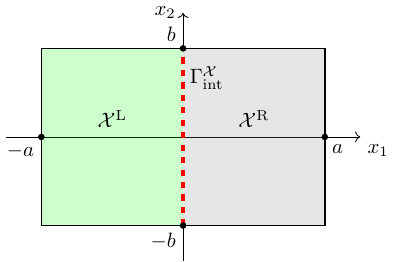}
        
    \caption{Domains $\calX^\text{L}$ and $\calX^\text{R}$.}
    \label{fig:X}
\end{figure}

We introduce the following setting: We divide $\calX$ (as in \cref{def:XY}) into
\[
\calX^\tL = (-a , 0) \times (-b , b), \quad  \calX^\tR = (0, a) \times (-b,b),
\]
and we define the interface segment $\Gamma^\calX_\I:=\{0\}\times(-b,b)$ (see \cref{fig:X}). The subscript \say{$\I$} should be read as \textit{interface}. 

We consider mappings $\m^\tL \in [C^1(\overbar{\calX^\tL})]^2$, $\m^\tR \in [C^1(\overbar{\calX^\tR})]^2$, such that,
\begin{subequations}
\label{eq:mix_mappings}
\begin{align}
    \label{subeq:mix_ellhyp}
    \det \D \m^\tL>0 \;\; \on \overbar{\calX^\tL}, 
    \quad
    &\det \D \m^\tR<0 \;\; \on \overbar{\calX^\tR}, 
    \\
    \label{subeq:conserv}
    \curl \m^\tL = 0, \text{ with  $\D\bmu^\tL$ SPD} \;\; \on \calX^\tL, 
    \quad 
    &\curl \m^\tR = 0 \;\;\on \calX^\tR,
    \\
    \label{subeq:continuity}
    \m^\tL = \m^\tR \;\; 
    \on &\GamintX .
\end{align}
\end{subequations}
By virtue of \cref{subeq:conserv,subeq:continuity} we can then construct $\m \in C(\calX)$ and $u \in C^1(\calX)$ via $\restr{\m}{\calX^\tL} = \m^\tL$, $\restr{\m}{\calX^\tR} = \m^\tR$ and $\nabla u = \m$. 

Ultimately, we want $\m$ to be a solution of \cref{eq:transport_problem} and in particular, $\m(\calX) = \calY$. However, it is important to note that the change in sign of $\det \D\m$ across $\GamintX$ is incompatible with global injectivity. Intuitively, $\m$ folds across the interface, with $\m(\GamintX)$ representing the resulting crease of the fold.

Motivated by this non-injectivity, we choose to enforce a global 2-to-1 structure on $\m$, that is, we assume both $\m^\tL: \overbar {\calX^\tL} \to \overbar \calY$ and $\m^\tR: \overbar {\calX^\tR} \to \overbar \calY$ to be diffeomorphisms, meaning that the preimage of each point in $\calY$ has precisely two points, one in $\calX^\tL$ and the other one in $\calX^\tR$. This choice comes particularly useful when searching for periodic mapping solutions as will be seen in Section \ref{sec:periodicity}.

\begin{figure}[ht]
    \centering
    \subfloat[Domain $\calX$.]{
        \includegraphics{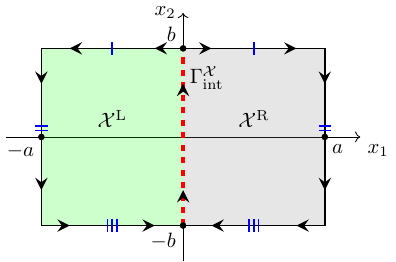}
    }
\hspace{5mm}
    \subfloat[Domain $\calY$.]{
        \includegraphics{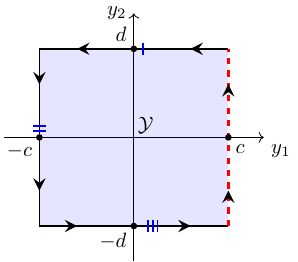}
    }
    \caption{Depiction of the boundary segmentation \cref{eq:boundary_segmentation_mix}}
    \label{fig:segmentation_mix}
\end{figure}

\begin{proposition}[Mixed Boundary Segmentation]\label{prop:boundary_segmentation_mix}
    Let $\m^\tL: \overbar {\calX^\tL} \to \overbar \calY$ and $\m^\tR: \overbar {\calX^\tR} \to \overbar \calY$ be diffeomorphisms satisfying \cref{eq:mix_mappings}. Then the following holds:
    \begin{subequations}
    \label{eq:boundary_segmentation_mix}
    \begin{alignat}{4}
    \label{eq:seg_mL}
    \m^\tL(\Gamma^{\calX^\tL}_\wrom) &= \Gamma^\calY_\wrom, \; \;
    \m^\tL(\Gamma^{\calX^\tL}_\nrom) &&= \Gamma^\calY_\nrom, \; \;
    \m^\tL(\Gamma^{\calX^\tL}_\erom) &&= \Gamma^\calY_\erom, \; \;
    \m^\tL(\Gamma^{\calX^\tL}_\srom) &&= \Gamma^\calY_\srom,
    \\
    \label{eq:seg_mR}
    \m^\tR(\Gamma^{\calX^\tR}_\wrom) &= \Gamma^\calY_\erom, \; \;
    \m^\tR(\Gamma^{\calX^\tR}_\nrom) &&= \Gamma^\calY_\nrom, \; \;
    \m^\tR(\Gamma^{\calX^\tR}_\erom) &&= \Gamma^\calY_\wrom, \; \;
    \m^\tR(\Gamma^{\calX^\tR}_\srom) &&= \Gamma^\calY_\srom.
    \end{alignat}
    \end{subequations}
\end{proposition}
\begin{proof}
    Since $\det \D\m^\tL > 0$ on $\overbar \calX^\tL$ and $\m^\tL=\nabla u^\tL$ for $u^\tL$ convex, then following the reasoning of Proposition \ref{prop:ell_seg} we obtain \cref{eq:seg_mL}. By Remark \ref{rem:boundary_orientation}, $\m^\tR$ inverts the orientation of curves $\boldsymbol\gamma :[0,1] \to \partial \calX$. Furthermore, since $\m^\tL = \m^\tR$ on $\GamintX$, then $\m^\tR(\Gamma^{\calX^\tR}_\wrom) = \m^\tL(\Gamma^{\calX^\tL}_\erom) = \Gamma^\calY_\erom$. By following cyclically along $\partial \calX^\tR$ and $\partial \calY$ we get,
    \[
    \m^\tR(\Gamma^{\calX^\tR}_\wrom) = \Gamma^\calY_\erom \quad \implies \quad
    \m^\tR(\Gamma^{\calX^\tR}_\nrom) = \Gamma^\calY_\nrom, \quad
    \m^\tR(\Gamma^{\calX^\tR}_\erom) = \Gamma^\calY_\wrom, \quad
    \m^\tR(\Gamma^{\calX^\tR}_\srom) = \Gamma^\calY_\srom.
    \]
    Indeed note that the image of traveling clockwise along $\partial\calX^\tR$ becomes anticlockwise along $\partial\calY$.
\end{proof}
The boundary segmentations given in Proposition \ref{prop:boundary_segmentation_mix} are depicted in \cref{fig:segmentation_mix}. The segments of $\calX^\tL$ and $\calX^\tR$ are marked in accordance with their image segment on $\partial \calY$. For example, both $\GamintX$ and $\Gamma^\calY_\erom$ are drawn with a red dashed line. Similarly, $\Gamma^{\calX^\tL}_\nrom$, $\Gamma^{\calX^\tR}_\nrom$ and their image $\Gamma^\calY_\nrom$ are all marked with one blue dash.
The arrows drawn in the figure highlight that orientation is preserved from $\calX^\tL$ to $\calY$  while from $\calX^\tR$ to $\calY$ orientation is reversed.

\paragraph{Transport problem.}
Let us recall that we want to find a solution $\m$ with this mixed elliptic-hyperbolic structure while solving the transport problem \cref{eq:transport_problem}, i.e., we want to transport the density $f$ on $\calX$ to the density $g$ on $\calY$, where $f,g$ satisfy \cref{eq:bounded_away,eq:total_energy_conservation}. This structure implies that $\m^\tL$ transports the density $\restr{f}{\calX^\tL}$ to some density on $\calY$, which we denote as $g^\tL$, while $\m^\tR$ transports the remaining energy $\restr{f}{\calX^\tR}$ to some other density $g^\tR$ on $\calY$. 
Informally, the idea is that these two density contributions on $\calY$ should add up to the desired density $g$. Mathematically we can write this as follows:
\begin{subequations}
\label{sys:forward_mix_MA}
    \begin{alignat}{2}
        \label{subeq:mix_MA}
        &\det \D \m^\tL = \frac{f^\tL}{g^\tL \circ \m^\tL} \,\, \on \calX^\tL, \quad 
        &&\det \D \m^\tR = -\frac{f^\tR}{g^\tR \circ \m^\tR} \,\, \on \calX^\tR,
        \\
        \label{subeq:mix_curl}
        &\curl \m^\tL = 0, \text{ with  $\D\bmu^\tL$ SPD }\; \on \calX^\tL, \quad 
        && \curl \m^\tR = 0 \;\; \on \calX^\tR,
        \\
        \label{subeq:mix_boundary}
        &\m^\tL(\calX^\tL) = \calY, && \m^\tR(\calX^\tR) = \calY,
    \end{alignat}
    and with the coupling conditions,
    \begin{equation}
        \label{subeq:mix_coupling}
        \m^\tL = \m^\tR \,\, \on \GamintX, \quad g^\tL + g^\tR = g, \,\text{ and }\; g^\tL, g^\tR>0,
    \end{equation}
\end{subequations}
where, 
\begin{equation}\label{eq:def_fLfR}
    f^\tL = \restr{f}{\calX^\tL}, \quad f^\tR = \restr{f}{\calX^\tR}.
\end{equation}
Equations \cref{subeq:mix_MA,subeq:mix_curl,subeq:mix_boundary} define two Monge--Ampère systems: an elliptic system for $\m^\tL$ and a hyperbolic system for $\m^\tR$. By \cref{subeq:mix_coupling}, these systems are coupled to form a mixed elliptic--hyperbolic system. The first coupling condition enforces continuity of the mapping across $\GamintX$, while the second ensures that the transported densities $g^\tL$ and $g^\tR$ combine to recover the target density $g$. 

\color{black}

With these conditions, one can verify (by applying a change of variables on $\calX^\tL$ and $\calX^\tR$) that the piecewise-defined mapping $\m \in [C(\calX)]^2$ given by $\restr{\m}{\calX^\tL} = \m^\tL$ and $\restr{\m}{\calX^\tR} = \m^\tR$ satisfies the transport equation \cref{eq:transport_problem}. Note also that, again by a change of variables, we have the following compatability condition
\begin{equation}\label{eq:subdomain_energyconservation}
    \int_{\calX^\tL} f^\tL\dx = \int_\calY g^\tL \dy, \quad \int_{\calX^\tR} f^\tR\dx = \int_\calY g^\tR \dy,
\end{equation}

\paragraph{Backward formulation.}
The additional unknowns $g^\tL$ and $g^\tR$ pose a big challenge in handling system \cref{sys:forward_mix_MA}. However, it turns out that this dependence can be removed by reformulating the system in terms of these inverse mappings. Indeed by Proposition \ref{prop:diffeo}, both $\m^\tL$ and $\m^\tR$ are invertible, so we can define
\begin{equation*}
    \bmu^\tL := (\m^\tL)^{-1}, \quad \bmu^\tR := (\m^\tR)^{-1}.
\end{equation*}
We then get the following system:
\begin{subequations}
\label{sys:backward_mix_MA}
    \begin{alignat}{2}
        \label{subeq:mix_MA_mu}
        &\det \D \bmu^\tL \, \left( f^\tL \circ \bmu^\tL \right) - \det \D \bmu^\tR \, \left( f^\tR \circ \bmu^\tR \right) &&= g  \;\;\;\; \on \calY,
        \\
        \label{subeq:nonvan_det}
        &\det \D \bmu^\tL > 0, && \det \D \bmu^\tR < 0,
        \\
        \label{subeq:mix_curl_mu}
        &\curl \bmu^\tL = 0, \text{ with  $\D\bmu^\tL$ SPD},\quad 
        && \curl \bmu^\tR = 0,
        \\
        \label{subeq:mix_boundary_mu}
        &\bmu^\tL(\calY) = \calX^\tL, && \bmu^\tR(\calY) = \calX^\tR,
    \end{alignat}
    and with the coupling conditions,
    \begin{equation}
        \label{subeq:mix_coupling_mu}
        \bmu^\tL = \bmu^\tR \,\, \on \Gamma^\calY_\erom, \quad \bmu^\tL(\Gamma^\calY_\erom) = \GamintX.
    \end{equation}
\end{subequations}
Equations \cref{subeq:mix_MA_mu,subeq:nonvan_det} are obtained through \cref{subeq:mix_MA}, the second coupling condition in \cref{subeq:mix_coupling}, and $\D \bmu^\ast(\y) = [\D \m^\ast (\bmu^\ast(\y))]^{-1}$ for $\ast \in \{\tL,\tR\}$. 

The global boundary conditions \cref{subeq:mix_boundary_mu} can be replaced by $\bmu^\tL(\partial \calY) = \partial \calX^\tL$ and $\bmu^\tR(\partial \calY) = \partial \calX^\tR$ due to Proposition \ref{prop:transport_boundary}.  Following the boundary segmentation \cref{eq:boundary_segmentation_mix} (see \cref{fig:segmentation_mix}) we can write these boundary conditions as Dirichlet boundary conditions  on the components of $\bmu^\tL=(\mu^\tL_1, \mu^\tL_2)$ and $\bmu^\tR=(\mu^\tR_1, \mu^\tR_2)$ as follows:
\begin{equation}\label{eq:dirichlet_mu}
    \begin{array}{ll}
       \mu_1^\tL = -a, \;\; \mu^\tR_1 = a  
       &\on \Gamma^{\calY}_\wrom,  
       \\
       \mu_2^\tL = -b = \mu^\tR_2   &\on \Gamma_\srom^{\calY},
    \end{array}
    \quad
    \begin{array}{ll}
        \mu_1^\tL = 0 = \mu^\tR_1  &\on \Gamma_\erom^{\calY}, 
        \\
        \mu_2^\tL = b = \mu^\tR_2   & \on \Gamma_\nrom^{\calY}.
    \end{array}
\end{equation}
These Dirichlet conditions are particularly useful in the numerical setting for ensuring that boundaries are mapped to boundaries.

Formulation \cref{sys:backward_mix_MA} is more advantageous compared to \cref{sys:forward_mix_MA} since it eliminates the unknowns $g^\tL$ and $g^\tR$, leaving only four unknown scalar functions $\bmu^\tL=(\mu^\tL_1, \mu^\tL_2)$ and $\bmu^\tR=(\mu^\tR_1, \mu^\tR_2)$. However, system \cref{sys:backward_mix_MA} is underdetermined, in the sense that we have four unknown scalar fields and only three first-order scalar partial differential equations, that is \cref{subeq:mix_MA_mu} and \cref{subeq:mix_curl_mu}.
This structural imbalance suggests non-uniqueness of solutions, which we illustrate in the following example.

\begin{example}[Non-Uniqueness]\label{ex:non_unique}
Consider system \cref{sys:backward_mix_MA} and let
\[
\calX = \calY = (-1,1)^2 , \quad \calX^\tL =(-1,0)\times(-1,1), \quad \calX^\tR = (0,1)\times (-1,1).
\]
Let also $f = g \equiv 1/4$. The idea is to construct solution pairs $(\bmu^\tL,\bmu^\tR)$ of the form
\[
\bmu^\tL:
\begin{cases}
    \calY  \to  \calX^\tL \\
     (y_1,y_2) \mapsto  \big(\alpha(y_1), y_2 \big)
\end{cases} \quad
\bmu^\tR:
\begin{cases}
    \calY  \to  \calX^\tR \\
     (y_1,y_2) \mapsto  (\beta(y_1),y_2)
\end{cases}
\]
with $\alpha:[-1,1] \to [-1,0]$, $\beta:[-1,1] \to [0,1]$  differentiable and surjective functions satisfying
\begin{equation}\label{eq:ac_conditions}
\alpha'>0, \;\; \alpha(-1) = -1, \;\; \alpha(1) = 0, \qquad \beta'<0, \;\; \beta(-1) = 1, \;\; \beta(1) = 0.
\end{equation}
This implies that
\cref{subeq:nonvan_det,subeq:mix_curl_mu,subeq:mix_boundary_mu,subeq:mix_coupling_mu} are satisfied. Finally \cref{subeq:mix_MA_mu} reduces to
\begin{equation}\label{eq:deriv_conditions}
    \alpha'(y_1) - \beta'(y_1) = 1.
\end{equation}
So by defining $\beta(y_1) := \alpha(y_1) - y_1 + 1$, we see that we can generate a different solution to \cref{sys:backward_mix_MA} for each $\alpha$ in the set
\[
\left\{
\alpha\in C^1([-1,1]): \quad \alpha(-1) = -1, \quad \alpha(1) = 0, \quad 0<\alpha'(y_1)<1, \quad y_1\in[-1,1]
\right\}.
\]
Indeed it is easy to check that \cref{eq:ac_conditions,eq:deriv_conditions} are fulfilled. 
\end{example}

\subsection{Additional constraint: periodicity of the mapping}\label{sec:periodicity}

Motivated by the underdetermined nature of \cref{sys:backward_mix_MA} and the interest in finding periodic reflectors (see \cref{fig:array_smooth}), we impose an additional constraint on the mapping $\m$: periodicity in the $x_1$-direction. This assumption allows $\m$ to be extended continuously across copies of $\calX$, and yields a smooth extension of the potential function $u$ along the $x_1$-axis. 

The periodicity condition is given by
\begin{equation}\label{eq:periodicity_mapping}
    \m^\tL(-a,x_2) = \m^\tR(a,x_2), \quad x_2\in (-b, b).
\end{equation}

To express this condition in the framework of the inverse mappings, we recall the boundary correspondence between $\calX$ and $\calY$ (see \cref{fig:segmentation_mix}). Under this identification, the periodicity condition is equivalent to enforcing equality of the second component of the inverse mappings on the corresponding boundary segment, namely
\begin{equation}\label{eq:periodic_mu}
    \mu_2^\tL = \mu_2^\tR \quad \on \Gamma^\calY_\wrom.
\end{equation}

Define the extended domain $\calX_N \subset \RI^2$ and the density $f_N: \calX_N \to (0,+\infty)$, for $N\in \mathbb N$, as
\begin{gather*}
    \calX_N := \big(-a,\, a+2Na\big)\times(-b,b), \\
    f_N(x_1 + 2ak, x_2) := f(x_1,x_2), \quad (x_1,x_2)\in \calX,\quad k=0,\ldots,N.
\end{gather*}
This corresponds to placing $\calX$ and $N$ more adjacent identical copies along the $x_1$-axis, each carrying the same density $f$.

Given a mapping $\m$ satisfying \cref{sys:forward_mix_MA,eq:periodicity_mapping}, the periodicity condition ensures that $\m$ extends continuously across these copies. We therefore define its periodic extension $\m_N$ on $\calX_N$ by
\[
\m_N(x_1 + 2ak, x_2) := \m(x_1,x_2), \quad (x_1,x_2)\in \calX,\quad k=0,\ldots,N.
\]
By construction, $\m_N$ is continuous on $\calX_N$. Moreover, for any $B \subset \calY$, we have
\[
\int_{\m_N^{-1}(B)} f_N \dx
= (N+1)\int_{\m^{-1}(B)} f \dx
= (N+1)\int_B g \dy,
\]
where the second equality follows from \cref{eq:transport_problem}. Hence, $\m_N$ transports the extended density $f_N$ onto a scalar multiple of $g$, preserving the relative spatial distribution while scaling the total mass.

It remains to extend the potential function. Since $\m = (m_1,m_2) = \nabla u$, integrating along vertical segments and using \cref{eq:periodicity_mapping} yields
\[
u(-a,x_2) - u(-a,-b)
= \int_{-b}^{x_2} m_2(-a,z)\,\mathrm{d}z
= \int_{-b}^{x_2} m_2(a,z)\,\mathrm{d}z
= u(a,x_2) - u(a,-b).
\]
Thus,
\begin{equation}\label{eq:deltau}
u(a,x_2) = u(-a,x_2) + \delta u, \quad x_2 \in (-b,b),
\end{equation}
where $\delta u := u(a,-b) - u(-a,-b)$.

We define the extension $u_N$ by
\[
u_N(x_1 + 2ak, x_2) := u(x_1,x_2) + k\,\delta u,
\quad (x_1,x_2)\in \calX,\quad k=0,\ldots,N.
\]
By construction, $\nabla u_N = \m_N$, and continuity of $\m_N$ implies $u_N \in C^1(\calX_N)$.

\begin{remark}[Optics interpretation]\label{rem:optic_insight_periodic}
The periodicity assumption is motivated by the optical interpretation of \cref{eq:transport_problem}. Interpreting $\calX$ as a light-emitting source, the domain $\calX_N$ corresponds to an array of identical sources arranged periodically in space. Such configurations can be used in illumination optics when the goal is to increase total light density while preserving the spatial structure of the illumination pattern.

From this perspective, constructing a periodic solution $\m$ allows the associated reflector potential $u$ to be computed only on a single cell $\calX$. The periodicity constraint then ensures that this local solution extends naturally and smoothly to the enlarged domain $\calX_N$, without introducing any additional compatibility conditions. In this way, the global reflector on $\calX_N$ is obtained from a single local computation, while remaining $C^1$ across all interfaces between adjacent cells (see \cref{fig:array_smooth}). The resulting transport map preserves the spatial distribution of the target density $g$, while scaling its density by a factor of $N+1$.
\end{remark}

\section{Least-squares method for the mixed MA system}\label{sec:mix_LSA}

We introduce an iterative least-squares scheme for solving system \cref{sys:backward_mix_MA}. Each iteration consists of three consecutive projection steps, where densities $g^\tL,g^\tR$, matrix functions $\PP^\tL,\PP^\tR$, and mappings $\bmu^\tL,\bmu^\tR$ are updated in sequence. The method is inspired by Prins~\cite{corien_paper}.

After obtaining sufficiently good approximations $(\bmu^\tL,\bmu^\tR)$ (which will be discussed in Section \ref{subsec:stopping_criteria}), we retrieve the forward mapping $\m:\calX \to \calY$  by inverting the approximations $\bmu\expseq{\tL}{n}$ and $\bmu\expseq{\tR}{n}$, i.e., we define $\m:\calX\to \RI^2$ as,
\begin{equation}\label{eq:recovery_mapping}
\restr{\m}{\calX^\tL}= \m^\tL, \quad \restr{\m}{\calX^\tR}= \m^\tR \quad \text{where } \,\, \m^\ast:= (\bmu\expseq{\ast}{n})^{-1},
\end{equation}
for $\ast \in \{\tL,\tR\}$.
It is then crucial to preserve injectivity across the iterations.

We then project $\m$ (via one final least-squares minimization) onto the space of conservative vector fields in order to obtain a potential function whose gradient is a good approximation of a solution to the transport problem \cref{eq:transport_problem}.

\subsection{Method description}

We start by choosing appropriate initial approximations $\bmu\expseq{\tL}{0},\bmu\expseq{\tR}{0}$ (which will be discussed at the end of this section) and at each iteration $n\in \mathbb N$, approximations $\bmu\expseq{\tL}{n+1}$ and $\bmu\expseq{\tR}{n+1}$ are computed from the previous iterate mappings via three consecutive least-squares minimization steps.  For compactness of notation we define the functional
\begin{equation*}
    J: [L^2(\calY) ]^{2\times 2}\times [L^2(\calY) ]^{2\times 2} \to \RI, \quad  J(\boldsymbol{D},\PP) := \frac{1}{2}\int_\calY | \boldsymbol{D} - \PP |_\text{F}^2\dy,
\end{equation*}
where $| \cdot |_\mathrm F$ is the Frobenius matrix norm given by $| \boldsymbol{A}|_\mathrm F=\sqrt{\text{tr}(\boldsymbol{A}^\mathrm{T} \boldsymbol{A}) }$ for $\boldsymbol{A}\in \RI^{2\times 2}$. 

We begin each iteration by computing the functions $h\expseq{\tL}{n},h\expseq{\tR}{n}$
\begin{equation}\label{eq:def_h1h2}
h\expseq{\tL}{n} := \big(f^\tL\circ\bmu\expseq{\tL}{n} \big)\,\det \D \bmu\expseq{\tL}{n}, \quad
h\expseq{\tR}{n} := -\big(f^\tR\circ\bmu\expseq{\tR}{n} \big)\,\det \D \bmu\expseq{\tR}{n},
\end{equation}
where $f^\tL,f^\tR$ are as in \cref{eq:def_fLfR} and $\bmu\expseq{\tL}{n},\bmu\expseq{\tR}{n}$ are the current iterate mappings.
The residual of \cref{subeq:mix_MA_mu} is then given by $h\expseq{\tL}{n}+h\expseq{\tR}{n}-g$ . The chosen notation $h\expseq{\ast}{n}$ (instead of $g\expseq{\ast}{n}$) is deliberate and highlights that $h\expseq{\tL}{n}+h\expseq{\tR}{n}$ may differ from $g$. 

The three minimization steps that comprise each iteration are the following:
\begin{subequations}
\label{eq:iterative_scheme}
\begin{align}
    \label{eq:minimization_g}
    (g\expseq{\tL}{n+1},g\expseq{\tR}{n+1}) &= \argmin_{(g^\tL,g^\tR)\in \calG_\varepsilon} \tfrac{1}{2}\|g^\tL-h\expseq{\tL}{n}\|^2_{L^2(\calY)} + \tfrac{1}{2}\|g^\tR-h\expseq{\tR}{n}\|^2_{L^2(\calY)},
    \\
    \label{eq:minimization_P}
    \PP^{\ast,(n+1)} &= \argmin_{\PP^\ast \in \calP\expseq{\ast}{n+1}}   J(\D \bmu\expseq{\ast}{n},\PP^\ast), \quad \ast \in \{\tL, \tR\},
    \\
    \label{eq:minimization_mu}
    (\bmu\expseq{\tL}{n+1},\bmu\expseq{\tR}{n+1})  &= \argmin_{(\bmu^\tL,\bmu^\tR) \in \calM } J\big(\D \bmu^\tL, \PP\expseq{\tL}{n+1} \big) + J\big(\D \bmu^\tR, \PP\expseq{\tR}{n+1} \big),
\end{align}
\end{subequations}
and where the admissible spaces $\calG_\varepsilon$, $ \calP\expseq{\ast}{n+1}$ are given by,
\begin{subequations}
\begin{align}
    \label{eq:set_h1h2}
     &\calG_\varepsilon = \left\{
    (g^\tL,g^\tR)\in [L^2(\calY)]^2 \,:
    \;\; \int_\calY g^\ast \dy = \int_{\calX^\ast} f^\ast\dx, \quad 
    \begin{array}{l}
         g^\ast \geq \varepsilon g,    \\
         g^\tL+g^\tR = g
    \end{array}
    \right\},
    \\
    \label{eq:LS_set_Past}
    &\calP\expseq{\ast}{n+1} = \left\{
         \PP \in [L^2(\calY)]^{2\times 2} : \;\; \det \PP = \eta^\ast \frac{g\expseq{\ast}{n+1}}{f^\ast\circ\bmu\expseq{\ast}{n}}, \; \PP = \PP^\mathrm{T} 
         \right\} ,
    \\
    \label{eq:def_M_mu}
    &\calM = \Big\{ (\bmu^\tL, \bmu^\tR) \in [H^1(\calY)]^{4}: \; (\bmu^\tL,\bmu^\tR) \text{ satisfy } \cref{eq:dirichlet_mu}, \; \bmu^\tL = \bmu^\tR \;\;\on \Gamma^\calY_\erom  \Big\},
\end{align}
\end{subequations}
for $\ast \in \{\tL,\tR\}$ where $\eta^\tL=1$ and $\eta^\tR=-1$. Here, $\varepsilon>0$ is a small parameter which bounds the densities $g^\tL$ and $g^\tR$ away from 0.

An intuition for the steps in \cref{eq:iterative_scheme} is the following. The minimization problem \cref{eq:minimization_g} projects $h\expseq{\tL}{n}$, $h\expseq{\tR}{n}$ onto the space of pairs of densities summing to $g$ (as desired) while also satisfying the integral constraint as in \cref{eq:subdomain_energyconservation}. In \cref{eq:minimization_P}  we project $\D \bmu\expseq{\ast}{n}$ onto $\calP\expseq{\ast}{n+1}$. The determinant constraint encodes the Monge--Ampère relation corresponding to the updated density $g\expseq{\ast}{n+1}$ while the symmetry constraint $\PP = \PP^\text{T}$ comes from wanting the mappings to be conservative. In \cref{eq:minimization_mu} we reconstruct mappings whose Jacobians best match $\PP\expseq{\ast}{n+1}$ while enforcing the transport boundary condition via \cref{eq:dirichlet_mu}, and the continuity condition $\bmu^\tL = \bmu^\tR$ on $\Gamma^\calY_\erom$. The purpose of the parameter $\varepsilon$ in \cref{eq:set_h1h2} is to bound $|\det \PP\expseq{\ast}{n+1}|$ uniformly away from zero, which we found to reduce the occurrence of degenerate mappings in the numerical experiments.

\begin{proposition}[Nonemptiness of $\calG_\varepsilon$]
    The set $\calG_\varepsilon$ as in \cref{eq:set_h1h2} is nonempty for $\varepsilon \leq \min(\alpha,1-\alpha)$, where $\alpha := \frac{\int_{\calX^\tL} f^\tL}{\int_\calY g}$.
\end{proposition}
\begin{proof}
    Note that \cref{eq:bounded_away,eq:total_energy_conservation} implies
    $
    0< \int_{\calX^\tL} f^\tL \dx < \int_\calY g \dy
    $,
    therefore $\alpha \in (0,1)$. Let $g^\tL = \alpha g$ and $g^\tR = g-g^\tL = (1-\alpha)g$. It is immediate that $g^\ast \geq \varepsilon g$ for $\ast = \{\tL,\tR\}$, and a simple computation yields
    \[
    \int_\calY g^\tL =  \alpha\int_\calY g \dy = \int_{\calX^\tL} f^\tL \dx, \quad 
    \int_\calY g^\tR =  \int_\calX f\dx - \int_{\calX^\tL} f^\tL = \int_{\calX^\tR} f^\tR,
    \]
    where on the first equality of the right side we used \cref{eq:total_energy_conservation} and the left side equation. Therefore $(g^\tL,g^\tR) \in \calG_\varepsilon$.
\end{proof}

We proceed with an individual analysis of each step in \cref{eq:iterative_scheme}.

\textbf{$(g^\tL,g^\tR)$-Projection step.} Consider problem \cref{eq:minimization_g}. We note that $g^\tR$ is completely determined by $g^\tL$, and so $\calG_\varepsilon$ can be represented through the following set
\begin{align}\label{eq:def:GLeps}
\calG^\tL_\varepsilon  = \left\{
g^\tL\in L^2(\calY) :
\;\; \int_\calY g^\tL \dy = \int_{\calX^\tL} f^\tL\dx, \quad \varepsilon g \leq g^\tL \leq (1-\varepsilon)g \;\; \on \calY
\right\},
\end{align}
By writing $g^\tR = g - g^\tL$, then the objective functional in \cref{eq:minimization_g} can be written as
\begin{align}
    \nonumber
    \tfrac{1}{2}|g^\tL-h\expseq{\tL}{n}|^2 + \tfrac{1}{2}|g^\tR-&h\expseq{\tR}{n}|^2 =  \tfrac{1}{2}|g^\tL-h\expseq{\tL}{n}|^2 + \tfrac{1}{2}|g-g^\tL-h\expseq{\tR}{n}|^2  
    \\
    \label{eq:completed_squares_gL}
    &= \left|g^\tL-\frac{g+h\expseq{\tL}{n}-h\expseq{\tR}{n}}{2}\right|^2 + \tfrac{1}{4}|g-h\expseq{\tL}{n}-h\expseq{\tR}{n}|^2.
\end{align}
And so, by removing the term independent of $g^\tL$, problem \cref{eq:minimization_g} is equivalent to
\begin{align}
\label{eq:minimization_gL}
g\expseq{\tL}{n+1} = \argmin_{g^\tL \in \calG^\tL_\varepsilon}  \tfrac{1}{2} \left\|g^\tL-h\right\|^2_{L^2(\calY)} \dy, \quad \text{where} \quad
h := \tfrac{1}{2}\big( g + h\expseq{\tL}{n}-h\expseq{\tR}{n}\big).
\end{align}
We then obtain the following characterization of the exact projection.
\begin{proposition}\label{prop:gL_projection}
    Let $\calG^{\tL}_\varepsilon$ be as in \cref{eq:def:GLeps} for $\varepsilon$ small enough, and $h\in L^2(\calY)$. Then
    \begin{equation}\label{eq:functional_gL}
    \argmin_{g^\tL \in \calG^\tL_\varepsilon} \tfrac{1}{2}\left\|g^\tL-h\right\|^2_{L^2(\calY)},
    \end{equation}
    has a unique minimizer $\hat g\in \calG^\tL_\varepsilon$ described by
    \begin{equation}\label{def:minimizer_g}
        \hat g = z(\hat \lambda), \quad \text{where }\;\; z(\lambda):= \max\Big(\varepsilon g, \, \min\big((1-\varepsilon)g,\, h-\lambda\big) \Big),
    \end{equation}
    and where $\hat \lambda\in \RI$ is chosen such that $\int_\calY z(\hat\lambda) \dy = \int_{\calX^\tL} f^\tL \dx$.
\end{proposition}
\begin{proof}
    The description of $\hat g$ \cref{def:minimizer_g} follows the logic of the Karush–Kuhn–Tucker conditions \cite{kuhntucker} where $\lambda$ is the Lagrangian multiplier enforcing $\int_\calY \hat g \dy = \int_{\calX^\tL} f^\tL\dx$.
    We prove the existence of $\hat\lambda\in \RI$ by showing that the function
    \begin{equation}\label{eq:lagrangemult_gL}
    I(\lambda) :=\int_\calY z(\lambda) \dy - \int_{\calX^\tL} f^\tL\dx,
    \end{equation}
    admits a root. Note that $\lim_{\lambda \to +\infty} z(\lambda) = \varepsilon g$ and $\lim_{\lambda \to -\infty} z(\lambda) = (1-\varepsilon)g$. Consequently, nonemptiness of $\calG^\tL_\varepsilon$ implies $\lim_{\lambda \to +\infty} I(\lambda) \leq 0$ and $\lim_{\lambda \to -\infty} I(\lambda) \geq 0$. Therefore by continuity and monotonicity of $I$, there exists $\hat \lambda$ such that $I(\hat \lambda) = 0$.
    We prove now that $\hat g$ as in \cref{def:minimizer_g} is indeed a minimizer. For any $v\in\calG^\tL_\varepsilon$,
    \[
    \|v-h\|_2^2 = \|v-\hat g\|_2^2 + 2\langle v-\hat g,\hat g-h\rangle + \|\hat g-h\|_2^2,
    \]
    where $\langle\cdot,\cdot\rangle$ is the usual $L^2$ inner product. 
    We prove that $\langle v-\hat g,\hat g-h\rangle \ge 0$ which then implies $\|v-h\|_2^2\ge \|\hat g-h\|_2^2$. Note that $v,h\in \calG_\varepsilon^\tL$ implies $\int_\calY v - \int_\calY \hat g = 0$, therefore
    \[
    \int_\calY (v-\hat g)(\hat g - h) \dy = \int_\calY (v-\hat g)(\hat g - (h-\lambda)) \dy.
    \]    
    We can partition the integration domain in three disjoint sets,
    \[
    A_-=\{h-\lambda\le \varepsilon g\},
    \quad
    A_0=\{\varepsilon g<h-\lambda<(1-\varepsilon)g\},
    \quad
    A_+=\{h-\lambda\ge (1-\varepsilon)g\}.
    \]
    On $A_0$ the integral is $0$ since $\hat g = h-\lambda$. On $A_-$ we have $\hat g = \varepsilon g$, and since $v \geq \varepsilon g$ we see that $(v-\hat g)\big(\hat g - (h-\lambda)\big) \geq 0$. A similar argument also leads to $(v-\hat g)\big(\hat g - (h-\lambda)\big) \geq 0$ on $A_+$, therefore $\langle v-\hat g,\hat g-h\rangle \ge 0$ easily follows.
    Uniqueness of the minimizer follows from the strict convexity of the objective functional \cref{eq:functional_gL}.
\end{proof}
The updated $g\expseq{\tL}{n+1}$ is then obtained by computing the exact projection \cref{def:minimizer_g}
and the second density is recovered via
\begin{equation}
    \label{eq:recovery_gR}
     g\expseq{\tR}{n+1} = g - g\expseq{\tL}{n+1}.
\end{equation}

\textbf{$\PP^\ast$-Projection step.}
We now consider the projection problem \cref{eq:minimization_P}. Since the objective functional does not depend on derivatives of $\PP^\ast$, the minimization is made pointwise where for each $\y\in \calY$ we solve a problem of the type
\begin{equation}\label{eq:pointwise_minimizationP}
\PP^\ast(\y)= \argmin_{\PP=\PP^\mathrm{T}\in\mathbb R^{2\times 2},\,\det\PP=C}
|\PP- \boldsymbol{D}|_{\mathrm F}^2,\;\; \text{ where }
C =\eta^\ast \frac{g\expseq{\ast}{n+1}(\y)}{f^\ast\circ\bmu\expseq{\ast}{n}(\y)}
\end{equation}
for $\ast\in\{\tL,\tR\}$, where $\boldsymbol{D} = \D\bmu\expseq{*}{n}(\y)$, and $\eta^\tL=1$ and $\eta^\tR=-1$. This problem has appeared previously in \cite{corien_paper,maikel}. We present here a compact characterization of the solution that applies to both the elliptic case ($C>0$) and the hyperbolic case ($C<0$).

Let $\boldsymbol{D}_s = (\boldsymbol{D}+\boldsymbol{D}^\mathrm{T})/ {2}$ be the symmetric part of $\boldsymbol{D}$. Since $\PP$ is constrained to be symmetric, a direct computation shows
\[
|\PP-\boldsymbol{D}|_{\mathrm F}^2
=
|\PP-\boldsymbol{D}_s|_{\mathrm F}^2
+
\frac{(d_{12}-d_{21})^2}{2}.
\]
The second term is independent of $\PP$, and therefore we may assume without loss of generality that $\boldsymbol{D}$ is symmetric.

Let $
\boldsymbol{D}
=
\boldsymbol{Q}
\diag(\lambda_1,\lambda_2)
\boldsymbol{Q}^\mathrm{T}$
be an eigendecomposition of $\boldsymbol{D}$ with $\boldsymbol{Q}$ orthogonal. By invariance of the Frobenius norm under orthogonal transformations
\[
|\PP-\boldsymbol{D}|_{\mathrm F}
=
|\widetilde{\PP}-\operatorname{diag}(\lambda_1,\lambda_2)|_{\mathrm F},
\qquad
\widetilde{\PP}:=\boldsymbol{Q}^\mathrm{T}\PP\boldsymbol{Q}.
\]
The following result shows that the minimizers share the eigenvectors of $\boldsymbol{D}$.

\begin{proposition}\label{prop:}
Let $\boldsymbol{Q}^\mathrm{T} \boldsymbol{D} \boldsymbol{Q} = \diag(\lambda_1,\lambda_2)$ be an eigendecomposition of $\boldsymbol{D}$ with $\boldsymbol{Q}$ orthogonal and $\lambda_1,\lambda_2 \in \RI$. Then any minimizer $\boldsymbol{P}$ of \cref{eq:pointwise_minimizationP} is also diagonal in the eigenbasis $\boldsymbol{Q}$, that is $\boldsymbol{Q}^\mathrm{T} \boldsymbol{P} \boldsymbol{Q} = \diag(\rho_1,\rho_2)$, for some $\rho_1,\rho_2 \in \RI$.
\end{proposition}
\begin{proof}
    Let $\boldsymbol{D}=\boldsymbol{Q}\diag(\lambda_1,\lambda_2)\boldsymbol{Q}^\mathrm{T}$ be an eigendecomposition of $\boldsymbol{D}$ with $\boldsymbol{Q}$ orthogonal, and define $\widetilde{\PP}:=\boldsymbol{Q}^\mathrm{T}\PP\boldsymbol{Q}$. Since $\widetilde \PP$ remains symmetric we can write 
    \[
    \widetilde\PP = R(\theta) \diag(\rho,C/\rho) R(-\theta),\quad  \text{where} \quad R(\theta) = 
    \begin{bmatrix}
        \cos\theta & -\sin\theta \\ \sin\theta & \cos \theta
    \end{bmatrix},
    \]
    for some $\theta\in \RI$. A direct computation shows
    \[
    \left|\widetilde \PP - \diag(\lambda_1,\lambda_2)\right|_{\text{F}}^2 = (\lambda_1-C/\rho)^2 + (\lambda_2-\rho)^2  - 2\cos^2\theta(\rho - C/\rho)(\lambda_1-\lambda_2).
    \]
    Minimizing over $\theta$ we conclude that either $\cos^2\theta = 0$ or $\cos^2 \theta = 1$. In either case, it implies that $\widetilde \PP$ is diagonal.
\end{proof}

Consequently, it suffices to determine the eigenvalues of $\PP$. Since we demand $\det \PP = C$, then we can write $\widetilde{\PP}=
\diag\left(\rho,\frac{C}{\rho}\right).$ 
The projection problem reduces to the one-dimensional minimization
\begin{equation}\label{eq:minimization_rho}
    \argmin_{\rho\in\mathbb R\setminus\{0\}}
    F(\rho),
    \;\; \text{ with } \;\;
    F(\rho)
    :=
    (\rho-\lambda_1)^2
    +
    \left(\frac{C}{\rho}-\lambda_2\right)^2.
\end{equation}

\begin{remark}[Geometric interpretation and generic uniqueness]\label{rem:geometry_hyperbola}
The minimization problem \eqref{eq:minimization_rho} can be interpreted as finding the closest point to $(\lambda_1,\lambda_2)$ on the hyperbola $\{(\rho, C/\rho)\}_{\rho \in \mathbb{R}\setminus\{0\}}$ (see \cref{fig:hyperbola_rho}). This projection is well-defined whenever the closest point is unique.

A geometric argument implies that non-uniqueness can only occur on the symmetry axes of the hyperbola, that is, for the symmetric configurations $\lambda_1 = \lambda_2$ or $\lambda_1 = -\lambda_2$.
In particular, uniqueness holds for almost every $(\lambda_1,\lambda_2) \in \RI^2
$, and therefore \cref{eq:pointwise_minimizationP} is well-defined for $\boldsymbol{D}\in \RI^{2\times 2}$ a.e.. 
\begin{figure}[ht]
    \centering
    \includegraphics[width=0.35\linewidth]{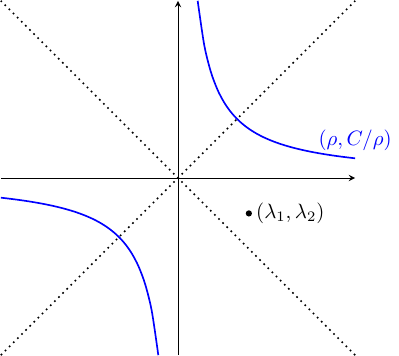}
    \caption{Geometric interpretation of \cref{eq:minimization_rho}.}
    \label{fig:hyperbola_rho}
\end{figure}
\end{remark}

We can compute the minimizers of $F$ via its stationary points. Differentiating yields
\begin{align}\label{eq:quartic_minimizationP}
F'(\rho) = 0 &\iff 2(\rho-\lambda_1) - 2\frac{C}{\rho^2}(C/\rho - \lambda_2) = 0
\iff \rho^4 - \lambda_1\rho^3 + C\lambda_2\rho - C^2 = 0.
\end{align}
The stationary points can therefore be computed analytically using, for example, Ferrari's method \cite[p.~32]{Tignol_galois}. Evaluating $F$ at all real stationary points yields a minimizer $\hat\rho$, from which the projection is reconstructed as
\begin{equation}
\label{eq:recovery_P}
\PP
=
\boldsymbol{Q}
\operatorname{diag}
\!\left(
\hat\rho,
\frac{C}{\hat\rho}
\right)
\boldsymbol{Q}^\mathrm{T}.
\end{equation}

\begin{remark}[SPD restriction of $\PP^\tL$]
    In the elliptic case ($\ast = \tL$), we can further restrict $\PP^\tL$ to be a symmetric positive definite (SPD) matrix. The reasoning behind it is that $\PP^\tL$ acts as an approximation of $\D\bmu^\tL$, and ultimately we want $\D\bmu^\tL$ to be pointwise an SPD matrix. 
    In order for $\PP^\tL$ to be SPD, the only adaptation that is needed is restricting $\rho$ to be positive in \cref{eq:minimization_rho}.
\end{remark}

\textbf{$\bmu^\ast$-Update step.} 
The final step of each iteration is \cref{eq:minimization_mu} which can be written as
\[
(\bmu\expseq{\tL}{n+1} , \bmu\expseq{\tR}{n+1}) = 
\argmin_{(\bmu^\tL,\bmu^\tR)\in \calM} J^{(n+1)}(\bmu^\tL,\bmu^\tR),
\]
where $\calM$ is defined in \cref{eq:def_M_mu} and $J^{(n+1)}:\calM \to \RI$ is defined as
\begin{equation}\label{def:Jn1}
J^{(n+1)}(\bmu^\tL,\bmu^\tR) :=
\frac{1}{2}\int_\calY \left( \left| \D\bmu^\tL - \PP\expseq{\tL}{n+1}\right|_F^2 + 
\left| \D\bmu^\tR - \PP\expseq{\tR}{n+1}\right|_F^2 \right) \dy.
\end{equation}

The minimizer is characterized by the vanishing of its first variation. We note that $\calM$ can be represented as 
\[
\calM = (\bmu^\tL,\bmu^\tR) + \calV,
\]
for some $(\bmu^\tL,\bmu^\tR) \in \calM$, where 
\begin{equation}\label{eq:def_calV}
\calV: = \left\{
(\pphi^\tL,\pphi^\tR)\in [H^1(\calY)]^{4}:
\begin{array}{ll}
   \varphi_1^\ast = 0
   & \on \Gamma^{\calY}_\wrom \cup \Gamma^\calY_\erom,  
   \\
   \varphi_2^\ast = 0   &\on \Gamma_\srom^{\calY} \cup \Gamma_\nrom^{\calY},
\end{array} 
\quad 
\pphi^\tL = \pphi^\tR \;\;\on \Gamma^\calY_\erom
\right\},
\end{equation}
for $\ast \in \{\tL , \tR \}$.
Through some simple computations we see that the first variation of $J^{(n+1)}(\bmu^\tL,\bmu^\tR)$ with respect to some $(\pphi^\tL,\pphi^\tR)\in \calV$ is given by
\begin{align}
\nonumber
\lim_{\varepsilon \to 0} \frac{1}{\varepsilon}&\left( {J^{(n+1)}\left(\bmu^\tL + \varepsilon\pphi^\tL,\bmu^\tR + \varepsilon\pphi^\tR\right) - J^{(n+1)}\left(\bmu^\tL,\bmu^\tR\right)} \right)
\\
\label{eq:weakform_minimizemu}
= &\int_\calY\Big[ (\D\bmu^\tL - \PP\expseq{\tL}{n+1}) \boldsymbol :\D\pphi^\tL +  (\D\bmu^\tR - 
\PP\expseq{\tR}{n+1}) \boldsymbol :\D\pphi^\tR   \Big]\dy = 0.
\end{align}
where $\boldsymbol{:}$ is the Frobenius inner product given by $\boldsymbol{A:B} ={\text{tr}(\boldsymbol{A}^\mathrm{T} \boldsymbol{B}) }$ for $\boldsymbol{A},\boldsymbol{B}\in \RI^{2\times 2}$.

\begin{proposition}\label{prop:uniqueness_muminimiz}
    There exists a unique minimizer $(\hat \bmu^\tL, \hat \bmu^\tR) \in \calM$ of $J^{(n+1)}$ (defined in \cref{def:Jn1}). Furthermore, \cref{eq:weakform_minimizemu} holds for every $(\pphi^\tL,\pphi^\tR) \in \calV$ if and only if $(\bmu^\tL,\bmu^\tR) = (\hat \bmu^\tL,\hat\bmu^\tR)$.
\end{proposition}
\begin{proof}
    Since $\calM$ is a Hilbert space and $J^{(n+1)}$ is strictly convex on $\calM$, it admits a unique minimizer. By the strict convexity we get additionally that this solution is the only stationary point of $J^{(n+1)}$. Therefore if \cref{eq:weakform_minimizemu} holds for some $(\bmu^\tL,\bmu^\tR) \in \calM$ and for every $(\pphi^\tL,\pphi^\tR) \in \calV$, then $(\bmu^\tL,\bmu^\tR)$ is the unique minimzer of $J^{n+1}$.
\end{proof}

Although \cref{eq:weakform_minimizemu} provides the variational
characterization of the minimizer, our numerical implementation will be based
on a finite-difference discretization of the corresponding strong form.
We therefore derive the Euler--Lagrange equations associated to
\cref{eq:weakform_minimizemu}.

Summing over the rows of the matrices in \cref{eq:weakform_minimizemu}, and through integration by parts and Gauss's theorem, we get
\begin{multline}\label{eq:mix_variational}
\sum_{i=1,2}\left[
    \oint_{\partial\calY} \varphi_i^\tL(\nabla \mu^\tL_i - \p_i^\tL)\bcdot \unitv{n}\ds \right.+ \oint_{\partial\calY} \varphi_i^\tR(\nabla \mu^\tR_i - \p_i^\tR)\bcdot \unitv{n}\ds
    \\
    \left.
    -\int_\calY (\Delta \mu^\tL_i - \nabla\cdot\p_i^\tL)\varphi_i^\tL\dy - \int_\calY (\Delta \mu^\tR_i - \nabla\cdot\p_i^\tR)\varphi_i^\tR\dy \right] = 0,
\end{multline}
where $\p^\ast_i$ is row $i$ of the matrix function $\PP^{\ast,(n+1)}$. Since the test functions are arbitrary in the interior of $\calY$, the Fundamental Lemma of Calculus of Variations implies, for $i=1,2$,
\[
\Delta \mu^\ast_i = \nabla\cdot\p_i^\ast \; \text{ on } \calY
\]

Let us look now at the boundary integrals. On the segments where $\varphi^\ast_i$ vanish we recover Dirichlet boundary conditions for $\mu^\ast_i$ (cf. \cref{eq:dirichlet_mu,eq:def_calV}).
From $\pphi^\tL=\pphi^\tR$ on $\Gamma^\calY_\erom$, the boundary integrals on $\Gamma^\calY_\erom$ in \cref{eq:mix_variational} from the L- and the R-part join as
\begin{equation}\label{eq:mix_coupled_neumann}
    \int_{\Gamma^\calY_\erom} \varphi_2^\tL\Big( \nabla \mu^\tL_2 - \p_2^\tL + \nabla \mu^\tR_2 - \p_2^\tR \Big)\bcdot \unitv{n} \ds.
\end{equation}
Since $\Gamma^\calY_\erom$ is the vertical east segment of $\partial \calY$, the outward normal is $\unitv{n}=(1,0)$. Recalling the continuity condition imposed on $(\bmu^\tL,\bmu^\tR)\in\calM$ in \cref{eq:def_M_mu}, we retrieve the following boundary conditions on $\Gamma^\calY_\erom$
\begin{equation}\label{eq:boundary_conditions_gamint}
\diff{\mu^\tL_2}{y_1} + \diff{\mu^\tR_2}{y_1} = p_{21}^\tL + p_{21}^\tR, \quad \mu^\tL_2 = \mu^\tR_2 \; \; \on \Gamma^\calY_\erom.
\end{equation}
On the remaining segments, since $\pphi^\tL$ and $\pphi^\tR$ are independent of each other we get standard Neumann boundary conditions.
\[
\diff{\mu^\ast_1}{y_2} = p^\ast_{12} \;\;\on \Gamma_\srom^{\calY} \cup \Gamma_\nrom^{\calY}, \quad \diff{\mu^\ast_2}{y_1} = p^\ast_{21} \;\; \on \Gamma^\calY_\wrom,
\]
for $\ast \in \{\tL,\tR\}$.
Combining all the equations, we obtain the following Poisson systems
for $(\bmu^\tL,\bmu^\tR)$,
\begin{subequations}\label{eq:poissons_mu}
    \begin{align}
        \Delta \mu^\tL_1 = \nabla \cdot \p_1^\tL, \quad 
        \Delta \mu^\tL_2 = \nabla \cdot \p_2^\tL, \quad
        \Delta \mu^\tR_1 = \nabla \cdot \p_1^\tR, \quad
        \Delta \mu^\tR_2 = \nabla \cdot \p_2^\tR, \;\;\on \calY,
        \\
        \label{eq:bc_poissonmu}
        \def\arraystretch{1.2}
        \begin{array}{llll}
             \mu^\tL_2=b=\mu^\tR_2, \;\;\; &  \displaystyle\partial_{y_2}{\mu^\tL_1} =p^\tL_{12}, \;\; & \displaystyle\partial_{y_2}{\mu^\tR_1} =p^\tR_{12}, \;\; &\on \Gamma^\calY_\nrom, 
            \\
            \mu^\tL_2=-b=\mu^\tR_2, \;\;\; & \displaystyle\partial_{y_2}{\mu^\tL_1} =p^\tL_{12}, \;\; &\displaystyle\partial_{y_2}{\mu^\tR_1} =p^\tR_{12}, \;\; &\on \Gamma^\calY_\srom,
            \\
            \mu^\tL_1=-a=-\mu^\tR_1, \; & \displaystyle\partial_{y_1}{\mu^\tL_2} =p^\tL_{21}, \;\; &\displaystyle\partial_{y_1}{\mu^\tR_2} =p^\tR_{21}, \;\; &\on \Gamma^\calY_\wrom, 
            \\
            \mu^\tL_1=0=\mu^\tR_1, \;\;\; & \mu^\tL_2=\mu^\tR_2, \;\;\; & \displaystyle \partial_{y_1}{\mu^\tL_2} + \partial_{y_1}\mu^\tR_2 = p_{21}^\tL+p_{21}^\tR,  &\on \Gamma^\calY_\erom, 
        \end{array}
    \end{align}
\end{subequations}
Note the coupling between $\mu^\tL_2$ and $\mu^\tR_2$ in the last line of \eqref{eq:bc_poissonmu} while there is no coupling condition between $\mu^\tL_1$ and $\mu^\tR_1$.

\begin{remark}\label{rem:minimization_periodic}
    If we also wish to enforce the periodicity of the mapping as suggested in Section \ref{sec:periodicity}, then this minimization step can be easily adapted. To add condition \cref{eq:periodic_mu}, we should alter the sets $\calM$ and ${\calV}$. The last condition of \cref{eq:def_M_mu} and \cref{eq:def_calV} should be changed respectively to 
    \[
    \mu^\tL_2 = \mu^\tR_2 \;\;\on \Gamma^\calY_\wrom \cup \Gamma^\calY_\erom,  \qquad 
    \varphi^\tL_2 = \varphi^\tR_2 \;\;\on \Gamma^\calY_\wrom \cup \Gamma^\calY_\erom.
    \]
    From here the variational techniques are identical to the ones just presented.
\end{remark}

\textbf{Summary of the least-squares method.}
We now sum up the iterative scheme introduced in this section. Let $(\bmu\expseq{\tL}{0},\bmu\expseq{\tR}{0})\in \calM$ be some initial approximations. Then the each iteration is given as follows:
\begin{itemize}
    \item[(i)] $g\expseq{\tL}{n+1},g\expseq{\tR}{n+1}$ are obtained via \cref{def:minimizer_g,eq:recovery_gR}.
    
    \item[(ii)] $\PP\expseq{\tL}{n+1},\PP\expseq{\tR}{n+1}$ are obtained by solving problem \cref{eq:minimization_P} pointwise. This is done by computing the zeros of \cref{eq:quartic_minimizationP} to find the minimizer $\hat \rho$ of \cref{eq:minimization_rho}. Then we recover $\PP\expseq{\ast}{n+1}(\y)$ via \cref{eq:recovery_P}.
    
    \item[(iii)] $\bmu\expseq{\tL}{n+1},\bmu\expseq{\tR}{n+1}$ can be obtained by solving the weak formulation \cref{eq:weakform_minimizemu}, although in our numerical implementation we solved \cref{eq:poissons_mu} via centered finite-differences.
\end{itemize}
Via Proposition \ref{prop:gL_projection}, step (i) is uniquely defined. Furthermore, from Remark \ref{rem:geometry_hyperbola} we see that step (ii) is uniquely defined if $|\lambda_1| \neq |\lambda_2|$ where $(\lambda_1,\lambda_2)$ are the eigenvalues of the symmetric part of $\D\bmu\expseq{\ast}{n}$, that is, if $(\lambda_1,\lambda_2)$ does not lie on the axis of the hyperbola (see \cref{fig:hyperbola_rho}). In our numerical experiments, the occurrence of these ill-posed settings was observed to be negligible. Step (iii) is also uniquely defined due to \cref{prop:uniqueness_muminimiz}.

\subsection{Recovery of the potential function}
After obtaining sufficiently good approximations $(\bmu^\tL,\bmu^\tR) = (\bmu\expseq{\tL}{n} , \bmu\expseq{\tR}{n})$ (see Section \ref{subsec:stopping_criteria}), we recover an approximation of the potential function $u$.

The recovery of the potential function requires the mappings
$\bmu^\tL,\bmu^\tR$ to be invertible.
We found that the parameter $\varepsilon>0$ in $\calG_\varepsilon$ (see \cref{eq:set_h1h2}) reduced the occurrence of non invertible mappings, however
we currently do not have a theoretical guarantee that the iterates
remain globally invertible. Therefore, in what follows, we assume that
the computed mappings are invertible.

We then define
\[
\m^\tL := (\bmu^\tL)^{-1}, \quad \m^\tR:= (\bmu^\tR)^{-1}.
\]
Since $\bmu^\tL=\bmu^\tR$ on $\Gamma^\calY_\erom$ and from the Dirichlet conditions \cref{eq:dirichlet_mu}, $\m^\tL$ and $\m^\tR$ agree
on the corresponding interface $\GamintX$, and we define the piecewise
mapping $\m$ by $\restr{\m}{\calX^\tL} = \m^\tL$ and $\restr{\m}{\calX^\tR} = \m^\tR$.

The potential function is obtained by solving the following least-squares minimization problem:
\[
u = \argmin_{\psi \in H^1(\calX)} I(\psi), \quad I(\psi) :=  \frac{1}{2} \int_{\calX^\tL} \big|\nabla\psi - \m^\tL \big|^2\dx + \frac{1}{2}\int_{\calX^\tR} \big|\nabla \psi - \m^\tR\big|^2\dx.
\]
Once more, using integration by parts, Gauss's theorem and the Fundamental Lemma of Calculus of Variations, we see that a stationary point of this functional satisfies
\begin{subequations}\label{eq:poisson_u}
\begin{alignat}{2}
    \Delta u &= \nabla\cdot \m^\tL \quad &&\on \calX^\tL,
    \\
    \Delta u &= \nabla\cdot \m^\tR \;\; &&\on \calX^\tR,
    \\
    \label{eq:poisson_u_neumannBC}
    \nabla u\cdot \unitv{n} &= \m \cdot \unitv{n} \;\; &&\on \partial \calX,
    \\
    \label{eq:poisson_u_new_interface}
    \diff{u}{x_1^-} &= \diff{u}{x_1^+} \quad &&\on \GamintX,
\end{alignat}
\end{subequations}
where, 
\begin{equation*}
    \diff{u}{x_1^\pm}(x_1,x_2) := \lim_{h \to 0^\pm} \frac{u(x_1+h,x_2) - u(x_1,x_2)}{h}.  
\end{equation*}
The interface condition in \cref{eq:poisson_u_new_interface}
arises from the cancellation of the interface terms in the first variation,
using the continuity of $\m$ across $\GamintX$. 

For the periodic case, a similar condition to \cref{eq:poisson_u_new_interface} can be added to allow for the smooth extension, that is,
\begin{equation}
    \diff{u}{x_1^+}(-a,x_2) = \diff{u}{x_1^-}(a,x_2), \quad \forall x_2\in (-b,b).
\end{equation}

\subsection{Stopping criteria}\label{subsec:stopping_criteria}

At each iteration, given the obtained approximate mappings $(\bmu\expseq{\tL}{n},\bmu\expseq{\tR}{n})$ we define the following three functions:
\begin{subequations}
\begin{alignat}{2}
    \label{eq:residual}
    &\epsilon_R^{(n)} := \Big|\,h\expseq{\tL}{n} + h\expseq{\tR}{n} - g \, \Big|,
    \\
    \label{eq:curl_residual}
    &\epsilon_\times^{(n)} := \left|\curl \bmu\expseq{\tL}{n} \right| + \left| \curl \bmu\expseq{\tR}{n}\right|,
    \\
    \label{eq:eps_mu}
    &\epsilon_\mu^{(n)}(\y) := \big\|\, \bmu\expseq{\tL}{n}(\y) - \bmu\expseq{\tL}{n-1}(\y) \, \big\|_{2} \; + \;\big\|\, \bmu\expseq{\tR}{n}(\y) - \bmu\expseq{\tR}{n-1}(\y) \, \big\|_2 ,
\end{alignat}
\end{subequations}
for $y\in \calY$, where $h\expseq{\tL}{n},h\expseq{\tR}{n}$ are as in \cref{eq:def_h1h2}. In \cref{eq:residual} we define the residual of \cref{subeq:mix_MA_mu} while in \cref{eq:curl_residual} we define $\epsilon_\times$ as the sum of the residuals of \cref{subeq:mix_curl_mu}. The third equation, \cref{eq:eps_mu}, simply represents the pointwise change in the mappings from iteration $n-1$ to iteration $n$.
To achieve global estimates we compute the $L^1$ norm of these functions, i.e.,
\begin{equation}\label{eq:def_taus}
\tau_R^{(n)} := \| \epsilon_R^{(n)} \|_{L^1}, 
\quad \tau_\times^{(n)} := \| \epsilon_\times^{(n)} \|_{L^1}, 
\quad \tau_\mu^{(n)} := \| \epsilon_\mu^{(n)} \|_{L^1} .
\end{equation}
If the algorithm converges to a discretized approximation, then the ideal choice of stopping criteria would be
\begin{equation}\label{eq:stop_crit_mu}
\tau_\mu^{(n)} < \delta_\mu,
\end{equation}
for a small threshold tolerance $\delta_\mu >0$. However, because we expect non-uniqueness of solutions it might also be unreasonable to expect convergence of the numerical approximations to a specific discretized solution. This will be highlighted in Example 1 in Section \ref{sec:results}. Moreover, due to non-uniqueness, even if we have an exact solution, we cannot expect the approximates to converge to this solution. Therefore, we monitor the quality of the numerical approximations through the behavior of the residuals $\tau_R,\tau_\times$.

We propose a more appropriate stopping criterion which is based on the relative change of the residual norms because once the residuals stabilize, then so does the accuracy of the approximations,
\begin{align}
    \label{eq:stop_crit}
    {|\tau_R^{(n-1)} - \tau_R^{(n)}|} <\delta_R^{\vphantom{n}}\tau_R^{(n)}, \quad
    {|\tau_\times^{(n-1)} - \tau_\times^{(n)}|} <\delta_\times^{\vphantom{n}}\tau_\times^{(n)},
\end{align}
for some chosen thresholds $\delta_R, \delta_\times> 0$ and where $n\in \mathbb{N}$ denotes the iteration number. 

\subsection{Numerical implementation} \label{subsec:implementation}
The derivations presented in the previous sections are carried out at
the continuous level. We now briefly describe the discretization used in
the numerical experiments.

We start by discretizing $\calY$ via a uniform $N_1\times N_2$ grid, and by choosing the initial approximations $(\bmu\expseq{\tL}{0},\bmu\expseq{\tR}{0})\in \calM$ (where $\calM$ is as in \cref{eq:def_M_mu}). For all numerical experiments presented in \Cref{sec:results} we used $\bmu\expseq{\tL}{0},\bmu\expseq{\tR}{0}$ to be the  the unique affine mappings satisfying the Dirichlet boundary conditions \cref{eq:dirichlet_mu}. More concretely,
\begin{equation}\label{eq:init_mappings}
\bmu^{\tL,(0)}(y_1,y_2) = \left(
(y_1-c)\frac{a}{2c} 
\;,\;
y_2\frac{b}{d}
\right),
\quad
\bmu^{\tR,(0)}(y_1,y_2) = \left(
(c-y_1)\frac{a}{2c}
\;,\;
y_2\frac{b}{d}
\right).
\end{equation}
Recall that $a,b,c,d>0$ are the values that define the rectangles $\calX$ and $\calY$ (see \cref{def:XY} and \cref{fig:segmentation_mix}).
However, we note that other initial mappings in $\calM$ can produce quite different numerical approximations which is related to system \cref{sys:backward_mix_MA} being underdetermined.
\color{black}

The $g^\tL$-projection step \cref{eq:minimization_g} is done by computing the zero $\hat \lambda$ of the function in \cref{eq:lagrangemult_gL}. The integral is approximated with the trapezoidal rule and $\hat \lambda$ is computed via the bisection method.

The $\PP^\ast$-projection step is done pointwise where the quartic equation in \cref{eq:quartic_minimizationP} is solved for each point in the discretization grid via the Ferrari method \cite[p.~32]{Tignol_galois}. 

The recovery of $\bmu^\ast$ is done by discretizing \cref{eq:poissons_mu} via second-order central finite differences. The resulting linear systems were solved in MATLAB using sparse LU decompositions stored as decomposition objects (because the matrices remain unchanged across the iterations).

The recovery of the potential function $u$ requires inverting $\bmu^\tL$ and $\bmu^\tR$ which is done via piecewise linear interpolation. System \cref{eq:poisson_u} is solved via second-order central finite differences on a uniform $2N_1\times N_2$ space discretization of $\calX$.

Lastly, the estimates $\epsilon_*,\tau_*$ in \cref{eq:residual,eq:def_taus}, for $*\in\{R,\times,\mu\}$, are also approximated via finite differences and the trapezoidal rule.

\section{Numerical results}\label{sec:results}

 In this section we test the performance of the algorithm proposed in \Cref{sec:mix_LSA} to solve system \cref{sys:backward_mix_MA} for two problem settings.
In both settings we take the sets $\calX = (-1,1) \times (-0.5,0.5)$ and $\calY= (-1,1) \times (-1,1)$, and the initial mappings  will be as in \cref{eq:init_mappings}.

\subsection{Example 1: Constant source density} For this first example we consider constant source densities $f
^\tL,f^\tR \equiv C_f$, and a simple target density $g(\y) = \frac{C_g}{|\y|^2 + 1}$ where $C_f,C_g\in \RI$ are such that
\begin{equation}\label{eq:numexamples_energies}
\int_{\calX^\ast} f^\ast \dx = \frac{1}{2}, \quad \int_\calY g \dy= 1.
\end{equation}
for $\ast\in\{\tL,\tR\}$.
We choose the clamping threshold value in \cref{def:minimizer_g} as $\varepsilon = 0.05$.

\begin{figure}[ht]
    \centering
    \hspace{-4mm}
    \subfloat[Results up to 1500 iterations.]
    {\label{fig:example1_taus1500_5050}
        \centering
        \includegraphics[height=0.3\textwidth]{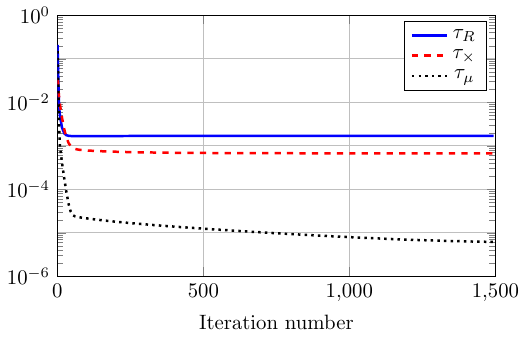}
    }
    \hspace{4mm}
    \subfloat[Results up to up to 25,000 iterations.]
    {\label{fig:example1_taus25000_5050}
        \centering
        \includegraphics[height=0.3\textwidth]{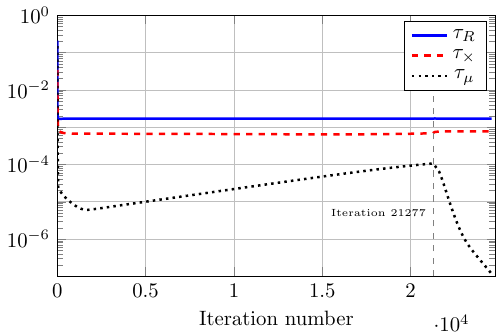}
    }
    \caption{Evolutions of $\tau_R,\tau_\mu,\tau_\times$ for a $50\times50$ discretization of $\calY$.}
    \label{fig:ex1_taus}
\end{figure}

We plot in \cref{fig:ex1_taus} the evolution of $\tau_R,\tau_\mu,\tau_\times$ for a uniform discretization grid of $50 \times 50$ points and for a fixed number of iterations. The first 1500 iterates are shown in \cref{fig:example1_taus1500_5050}, while \cref{fig:example1_taus25000_5050} shows up to iteration 25,000.

Initially, all three quantities decrease. Eventually, $\tau_R$ and $\tau_\times$ stabilize while $\tau_\mu$ increases gradually. This behavior persists for several thousand iterations, suggesting that the iterates continue to evolve while maintaning a comparable residual quality. Only after iteration 21,277 does $\tau_\mu$ reverse its trend and converge towards zero.

\begin{figure}[ht]
    \centering
    \subfloat[Grid lines on $\calY$.]
    {\label{fig:ex1_target}
        \centering
        \includegraphics[height=0.25\textwidth]{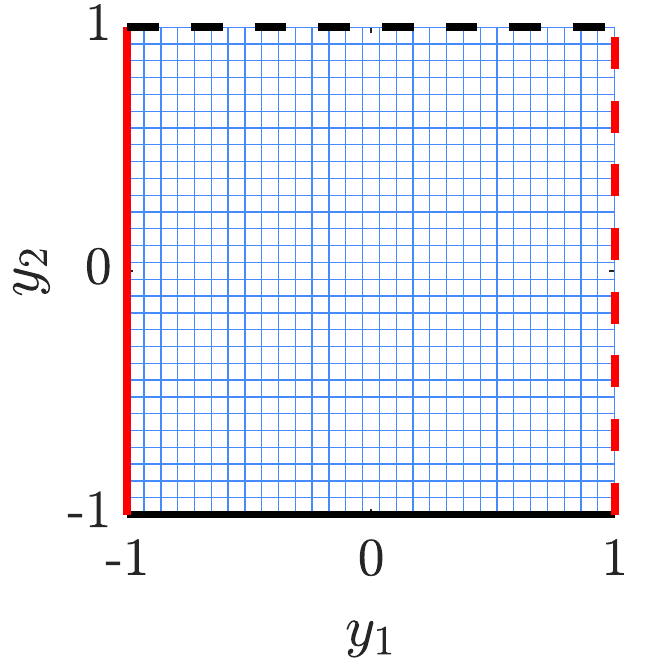}
    }

    \subfloat[$\bmu^\tL$ and $\bmu^\tR$ after 1500 iterations.]
    {\label{subfig:ex1_map1500}
        \centering
        \begin{tikzpicture}[spy using outlines={circle,black,magnification=5,size=1.5cm, connect spies}]
            \node {\includegraphics[width=0.40\textwidth]{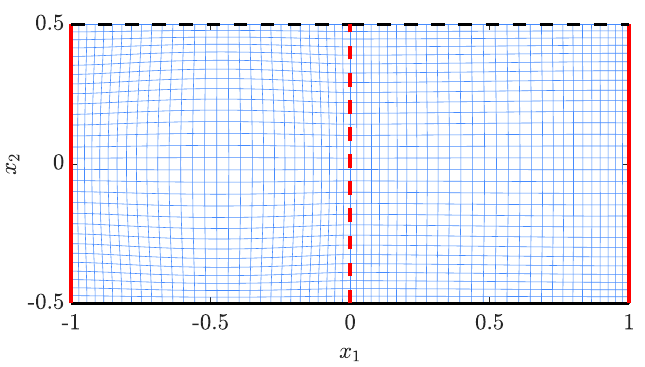}};
            \spy on (0.2,0.2) in node [left] at (2,0.3);
        \end{tikzpicture}
    }
    \subfloat[$\bmu^\tL$ and $\bmu^\tR$ after 21,277 iterations.]
    {\label{subfig:ex1_map21277}
        \centering
        \begin{tikzpicture}[spy using outlines={circle,black,magnification=5,size=1.5cm, connect spies}]
            \node {\includegraphics[width=0.40\textwidth]{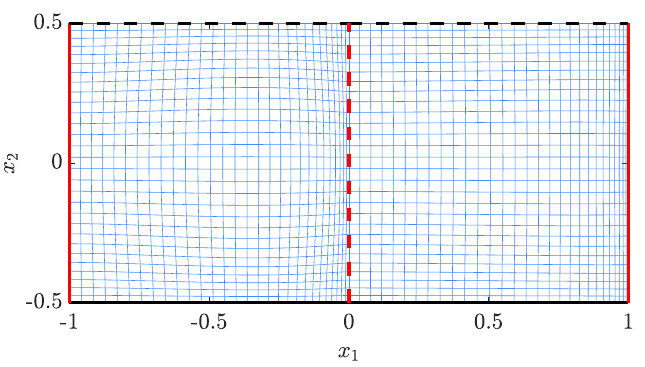}};
            \spy on (0.2,0.2) in node [left] at (2,0.3);
        \end{tikzpicture}
    }
    \caption{Inverse mappings $(\bmu^\tL,\bmu^\tR)$ on iterations 1500 and 21,277 for a $50\times50$ discretization.}
    \label{fig:ex1_mappings}

\end{figure}

The evolution of the inverse mappings provides insight into this phenomenon. \cref{fig:ex1_target} shows the target domain $\calY$ with equally spaced perpendicular grid lines, while \cref{subfig:ex1_map1500,subfig:ex1_map21277} display the images of these grid lines by $\bmu^\tL,\bmu^\tR$ at different iterations. The boundary segmentations are indicated by the markings on the boundaries of the domains.

Comparing \cref{subfig:ex1_map1500,subfig:ex1_map21277}, we observe a progressive compression of the vertical grid lines to the left of $x_1=0$ and to the left of $x_1=1$. This compression develops while $\tau_R$ and $\tau_\times$ remain nearly unchanged, indicating that the iterates drift along a family of approximate solutions rather than converging directly to a fixed point.

As a consequence of this drift, $|\det \D\bmu^\tL|$ decreases near $y_1=1$, while $|\det \D\bmu^\tR|$ decreases near $y_1=-1$. When these determinants become sufficiently small, the $\varepsilon$-clamping mechanism in the $g^\tL$-projection step becomes active, meaning that $h - \lambda$ falls outside the box constraints (cf. \cref{def:minimizer_g}). In this experiment, this first occurs at iteration 21,277 and persists throughout the simulation. This clamping coincides with the disappearance of the drifting behavior and the subsequent convergence of $\tau_\mu$ to zero.


\newpage

This example highlights the importance of the choice of stopping criterion. As observed previously, the quantity $\tau_\mu$ may continue to evolve long after $\tau_R$ and $\tau_\times$ have stabilized. Consequently, requiring the convergence of all three indicators can lead to a substantial increase in computational cost without noticeably improving the satisfaction of the transport and curl conditions. For this reason, we consider the stopping criterion \cref{eq:stop_crit}, based solely on $\tau_R$ and $\tau_\times$.

\begin{figure}[htbp]
    \centering
    \subfloat[{Convergence analysis of $\tau_R$ and $\tau_\times$ without periodicity condition.}]{
    \includegraphics[width=0.40\linewidth]{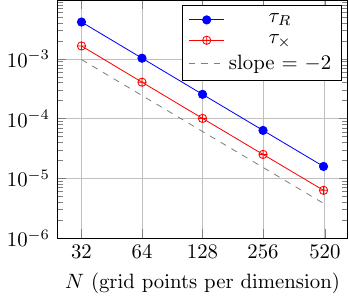}
    \label{fig:example1_convergence}
    }
    \hspace{10mm}
    \subfloat[Convergence analysis of $\tau_R$ and $\tau_\times$ with periodicity condition.]{
    \includegraphics[width=0.40\linewidth]{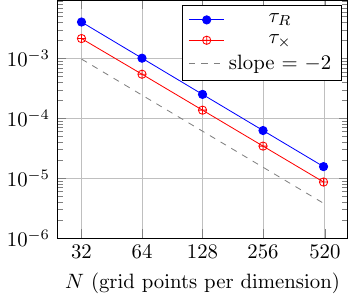}
    \label{fig:example1_convergence_periodic}
    }
    
    \subfloat[Potential function $u$ obtained from $\bmu^\tL,\bmu^\tR$ on a $512\times 512$ grid discretization of $\calY$ without periodicity condition]
    {\label{subfig:ex1_reflector}
        \centering
        \includegraphics[width=0.4\textwidth]{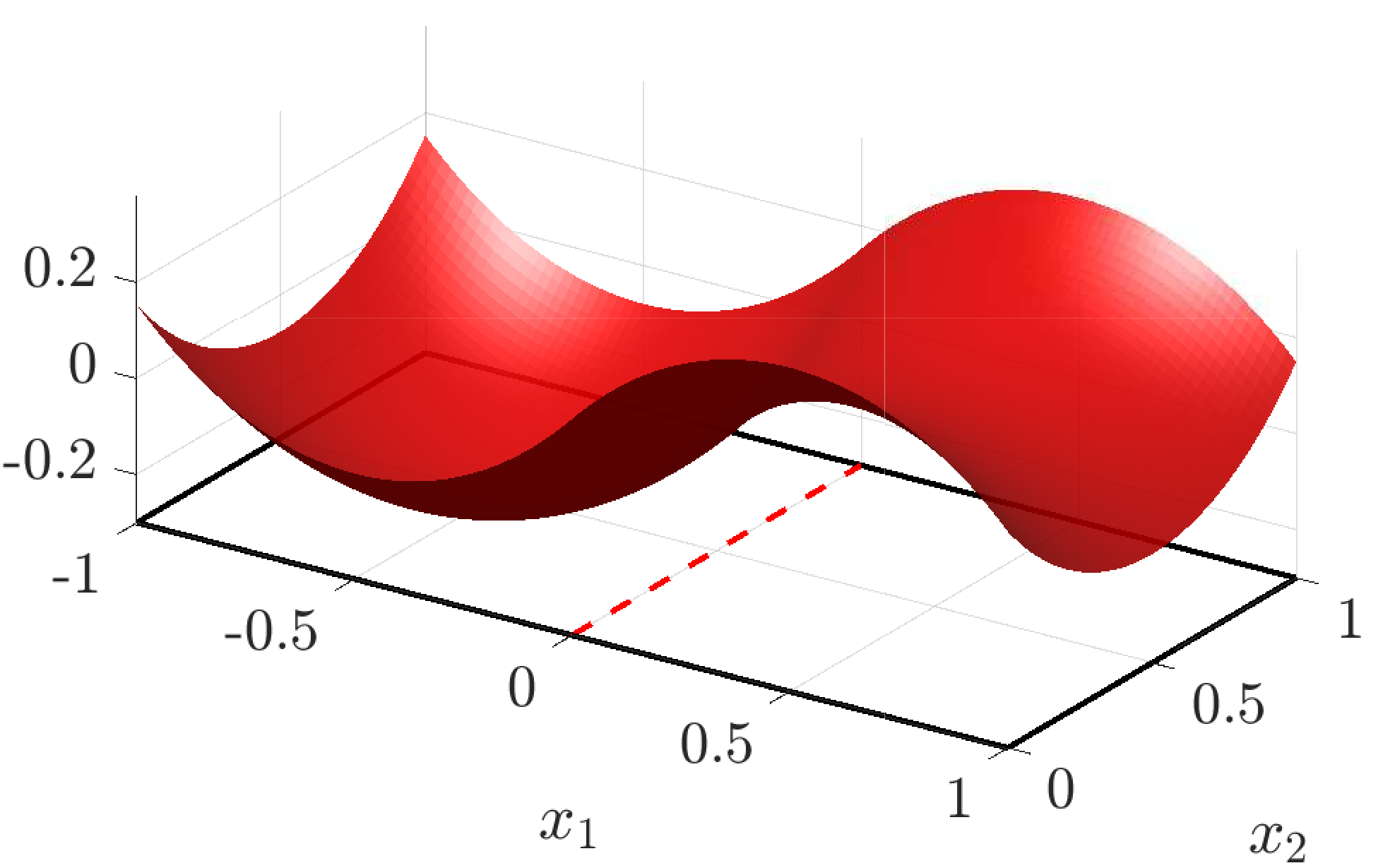}
    }
    \hspace{2mm}
    \subfloat[Extension of the potential function obtained from the periodicity condition.]
    {\label{subfig:ex1_periodic_reflector}
        \centering
        \includegraphics[width=0.43\textwidth]{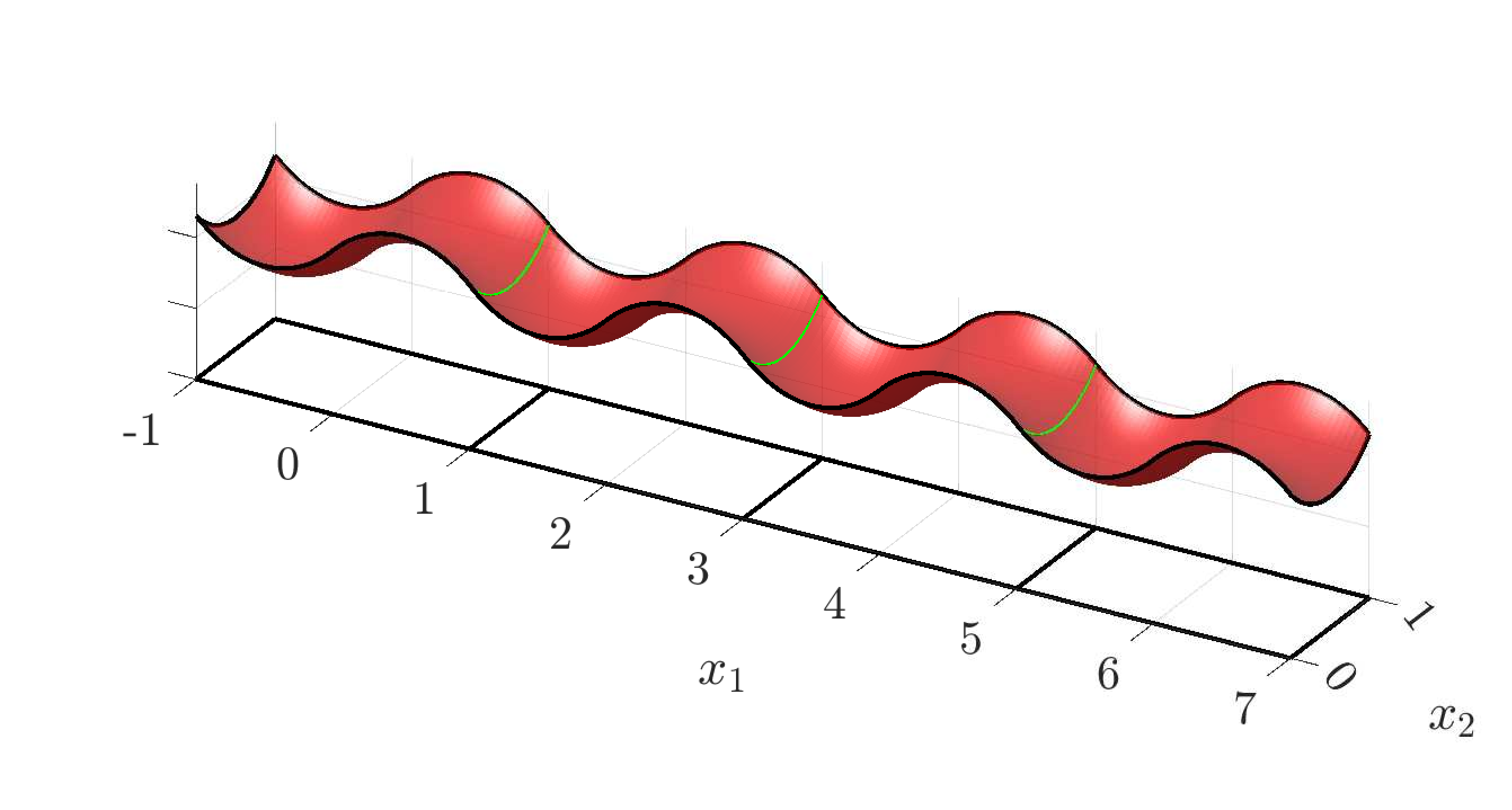}
    }
    
    \subfloat[Heatmap of $g^\calB$ (see \cref{eq:g_bins}).]
    {\label{subfig:density_g_ex1}
        \centering
        \includegraphics[width=0.4\textwidth]{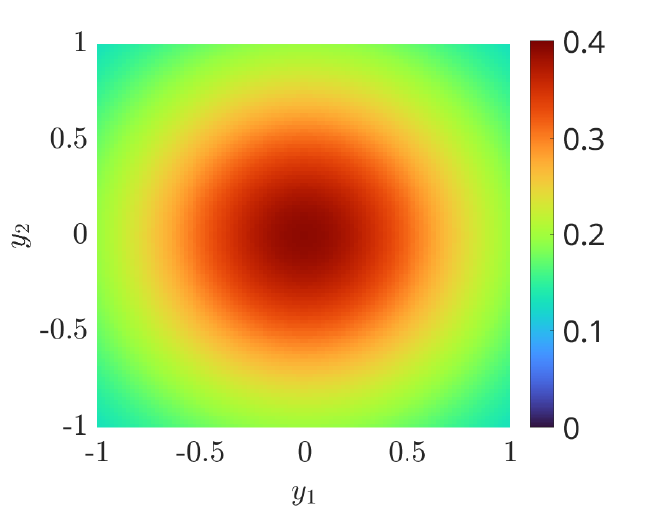}
    }
    \hspace{2mm}
    \subfloat[Heatmap of $|H_{N_s} - g^\calB|$ (see \cref{eq:bincounts,eq:g_bins}).]
    {\label{subfig:heatmap_diff_ex1}
        \centering
        \includegraphics[width=0.4\textwidth]{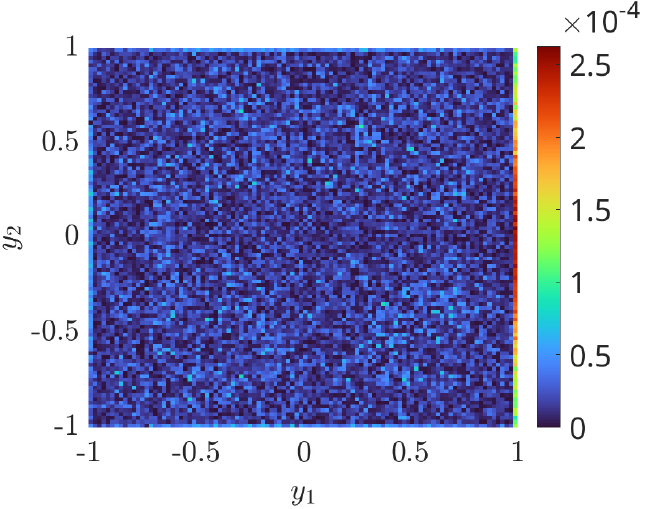}
    }
    \caption{Numerical results for Example 1.}
    \label{fig:heatmaps_ex1}
\end{figure}

\textbf{Convergence plots.}
So far we only analyzed the simulation of the algorithm for a fixed $50\times50$ grid discretization. Since mappings on each iteration are obtained by solving \cref{eq:poissons_mu} using second-order finite differences, the residuals $\tau_R$ and $\tau_\times$ are expected to depend on the grid discretization.

To investigate this dependence, we performed a convergence study by running the algorithm on uniform $N\times N$ grids with $N \in \{32,64,128,256,512\}$ using the stopping criteria \cref{eq:stop_crit} with relative tolerance values $\delta_R, \delta_\times=10^{-5}$. 

The results are displayed in \cref{fig:example1_convergence,fig:example1_convergence_periodic} where in the latter, we enforced the periodicity condition as discussed in Remark \ref{sec:periodicity} (cf. \cref{eq:periodic_mu}). In both cases we observe a second-order convergence of both residuals.

\textbf{Recovery of the potential function.}
The potential function $u$ obtained from the simulation of the algorithm with the $512\times512$ grid discretization without the periodicity condition is displayed in \cref{subfig:ex1_reflector}. As desired, the recovered potential is convex on $\calX^\tL$ and saddle-shape on $\calX^\tR$. 

In \cref{subfig:ex1_periodic_reflector}, we see the potential function and its extension to the domain $(-1,7)\times(0,1)$, obtained from the simulation with the periodicity condition. The figure illustrates that the reflector computed on the original domain extends smoothly across the replicated copies of $\calX$.

Consequently, the reconstructed gradient $\nabla u$ maps each boundary segment of the outer boundary $\partial \calX$ to the prescribed segment of $\partial \calY$, thanks to the boundary segmentation and the Neumann boundary condition \cref{eq:poisson_u_neumannBC}. On the interface $\GamintX$, however, since we impose the interface differentiability condition \cref{eq:poisson_u_new_interface} we don't have the additional freedom to additionally enforce $\partial_{x_1} u = c$ (where $c$ is as in \cref{def:XY}, see also \cref{fig:segmentation_mix}). The accuracy of this relation is therefore governed by the discretization errors and by the residual of the curl constraint.

\textbf{Sampling analysis.}
We assess how accurately the gradient of the numerical potential $u$ transports the source density $f$ onto the target density $g$ by a Monte Carlo sampling experiment.

We partition the target domain $\calY$ in a uniform grid of $M\in \mathbb N$ bins, $\mathcal{B}=\{B_k\}_{k=1}^M$. Then we draw $N\in \mathbb N$ samples uniformly distributed in $\calX$ (since the source density $f$ is constant) and compute their image via $\nabla u$. Subsequently we define the empirical histogram density
\begin{equation}\label{eq:bincounts}
    H_{N_s}(\y) := \displaystyle\frac{C_{k}}{N|B_{k}|}, \quad \text{if }\y \in B_{k},
\end{equation}
where $C_k$ denotes the number of mapped samples falling into $B_k$. We compare this with the bin-averaged target density
\begin{equation}\label{eq:g_bins}
    g^\calB(\y) :=\frac{1}{|B_{k}|} \int_{B_k} g\dy, \quad \text{if }\y \in B_{k}.
\end{equation}
In \cref{subfig:heatmap_diff_ex1} we see the results of this experiment with a $100\times100$ binning and $N_s=10^9$ samples, using the potential displayed in \cref{subfig:ex1_reflector}. The results show a small overall discrepancy between $H_{N_s}$ and $g^\mathcal{B}$, with a narrow localized increase in error near $\Gamma^\calY_\mathrm{int}=\Gamma^\calY_\erom$. This is consistent with the interface effect discussed previously, where $\partial_{x_1}u$ does not exactly match the prescribed boundary value. 
We omit here the results from the sampling analysis for the reflector obtained in the periodic setting as it was indistinguishable from \cref{subfig:heatmap_diff_ex1}.

\subsection{Example 2: Discontinuous target density} 

We consider again constant source densities $f
^\tL,f^\tR \equiv C_f$ and a discontinuous target density given by
\[
g(y_1,y_2) = \begin{cases}
    C_g(y_1^2 + 1) \quad &\text{if } (y_1-1)^2 + (y_2 + 0.6)^2 \geq 1,
    \\
    C_g/2          \quad &\text{if }  (y_1-1)^2 + (y_2 + 0.6)^2 < 1,
\end{cases}
\]
where the discontinuity is located at the to a arc of radius 1 centered at $(1,-0.6)$, and $C_f$ and $C_g$ are such that \cref{eq:numexamples_energies} holds. The discontinuity introduces a mismatch with the uniform Cartesian discretization and provides a more challenging test for the scheme.
We set the clamping threshold as $\varepsilon = 0.3$ and use stopping tolerances $\delta_R,\delta_\times = 10^{-5}$. 

\begin{figure}[htbp]
    \centering
    \subfloat[$(\bmu^\tL,\bmu^\tR)$ obtained without periodicity condition.]
    {\label{subfig:exdisc_mu}
        \centering
        \includegraphics[width=0.45\textwidth]{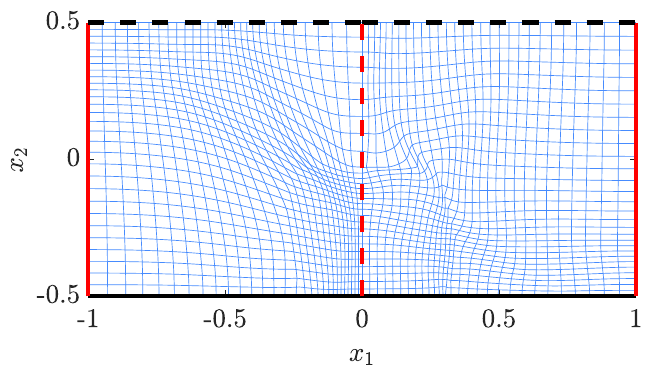}
    }
    \hspace{3mm}
    \subfloat[$(\bmu^\tL,\bmu^\tR)$ obtained with periodicity condition]
    {\label{subfig:exdisc_mu_periodic}
        \centering
        \includegraphics[width=0.45\textwidth]{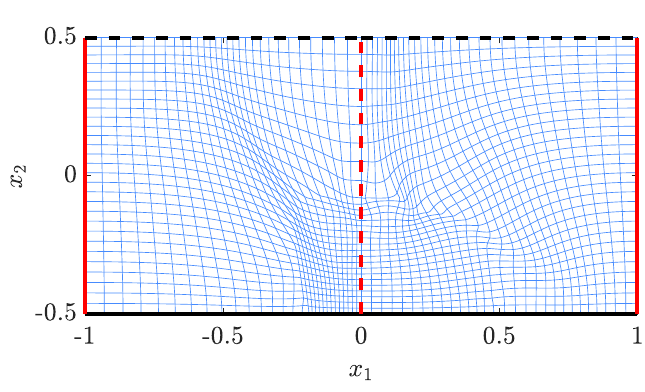}
    }
    \caption{Inverse mappings $(\bmu^\tL,\bmu^\tR)$ for Example 2 for a $512\times512$ discretization.}
    \label{fig:mappings_ex2}
\end{figure}

The mappings $(\bmu^\tL,\bmu^\tR)$ from the simulation with a $512\times 512$ spatial discretization, with and without the periodicity condition are shown in the first row of \cref{fig:mappings_ex2}. Note that in \cref{subfig:exdisc_mu_periodic}, each horizontal blue line intersects $\{-1\} \times [-0.5,0.5] = \Gamma^\calX_\wrom$ and $\{1\} \times [-0.5,0.5] = \Gamma^\calX_\erom$ at points with precisely the same $x_2$ value, which is a direct consequence of the periodicity condition. Indeed, \cref{subfig:exdisc_mu} differs clearly from \cref{subfig:exdisc_mu_periodic} in this aspect.
\color{black} 
It is also important to observe the loss of regularity of the mappings in both figures.
Due to the discontinuity of $g$, the Jacobians $\D \bmu^\tL$ and $\D \bmu^\tR$ are expected to lose regularity near the discontinuity region, while the mappings themselves remain continuous. This loss of smoothness is clearly visible in the deformation patterns in \cref{fig:mappings_ex2}.

\begin{figure}[htbp]
    \centering
    \subfloat[Convergence analysis of $\tau_R$ and $\tau_\times$ without periodicity condition.]
    {\label{subfig:exdisc_conv_non_periodic}
        \centering
        \includegraphics[width=0.4\textwidth]{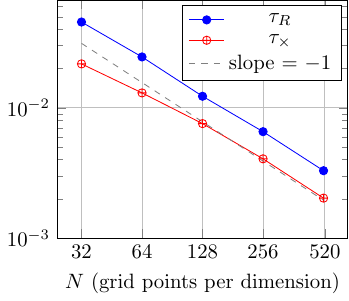}
    }
    \hspace{3mm}
    \subfloat[Convergence analysis of $\tau_R$ and $\tau_\times$ with periodicity condition.]
    {\label{subfig:exdisc_conv_periodic}
        \centering
        \includegraphics[width=0.4\textwidth]{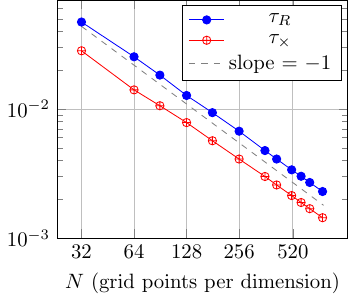}
    }
    
    \subfloat[Potential function $u$ obtained from $\bmu^\tL,\bmu^\tR$ without periodicity condition.]
    {\label{subfig:exdisc_potential}
        \centering
        \includegraphics[width=0.40\textwidth]{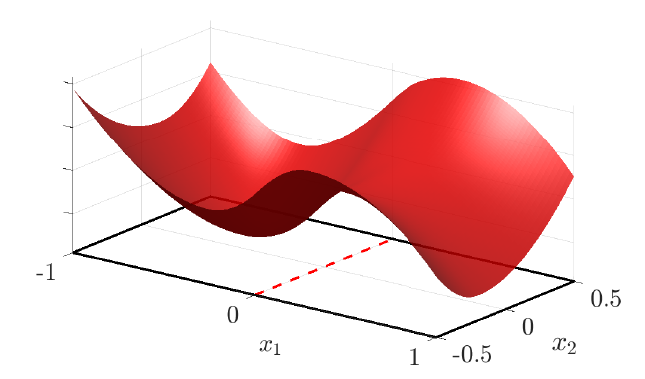}
    }
    \hspace{3mm}
    \subfloat[Extension of the potential function obtained from the periodicity condition.]
    {\label{subfig:exdisc_potential_periodic}
        \centering
        \includegraphics[width=0.40\textwidth]{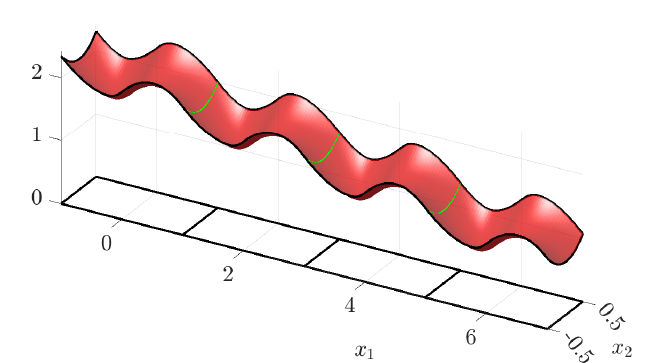}
    }
    
    \subfloat[Heatmap of $g^\calB$ (see \cref{eq:g_bins}).]
    {\label{subfig:exdisc_heatmap_g}
        \centering
        \includegraphics[width=0.43\textwidth]{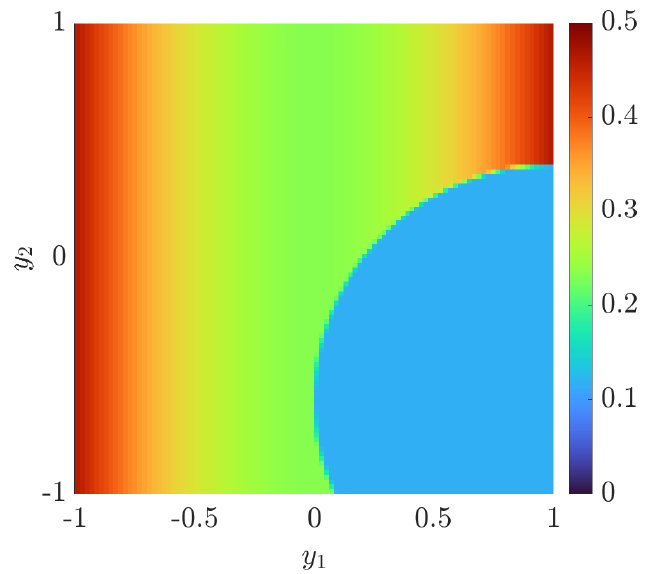}
    }
    \subfloat[Heatmap of $|H_{N_s} -g^\calB|$(see \cref{eq:bincounts,eq:g_bins}).]
    {\label{subfig:heatmap_diff_exdisc}
        \centering
        \includegraphics[width=0.43\textwidth]{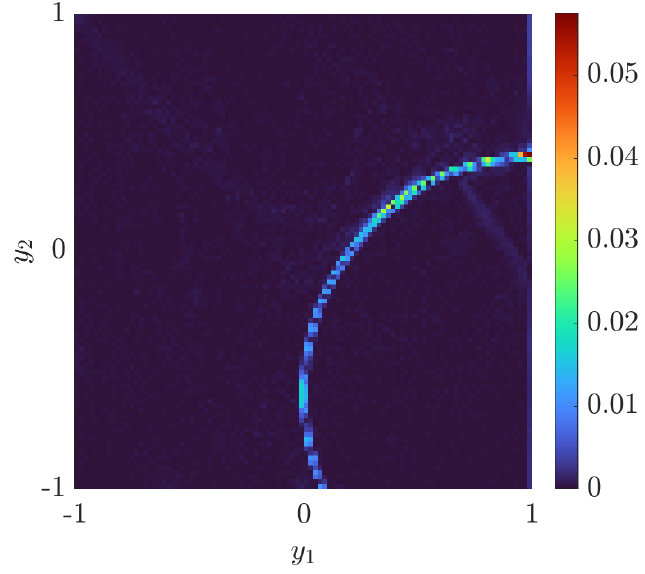}
    }
    \caption{Numerical Results for Example 2.}
    \label{fig:results_ex2}
\end{figure}

\textbf{Convergence analysis.}
In \cref{subfig:exdisc_conv_non_periodic,subfig:exdisc_conv_periodic} we display the results from the convergence study with $N \in \{32,64,128,256,512\}$. These results show first-order convergence of $\tau_R$ and $\tau_\times$, which is consistent with the presence of discontinuities in $g$. Indeed, these error estimates are computed via second-order finite-differences, however, a solution $\bmu^\tL,\bmu^\tR$ is no longer expected to be smooth in the discontinuity region of $g$. Therefore, near this region the pointwise discretization errors are large since they fail to capture the sudden jump of $g$. When integrating over the whole domain, these errors contribute as $\mathcal{O}(h)$ to the overall global residual estimation, where $h \propto 1/N$ is the grid spacing.

\textbf{Recovery of the potential function.}
\cref{subfig:exdisc_potential,subfig:exdisc_potential_periodic} display the potential function $u$ obtained from the simulation of the algorithm with the $512\times512$ grid discretization with and without the periodicity condition. However, we note that the $\delta u$ obtained (see \cref{eq:deltau} in Section \ref{sec:periodicity}) was nonzero, therefore the periodic extension is slanted. In this present work we were not able to control this parameter.

\textbf{Sampling analysis.}
We perform the same Monte Carlo analysis as in Example~1 using a $100\times100$ binning and $N_s=10^9$ samples, based on the mappings computed on a $512\times512$ grid without periodicity. The histogram $H_{N_s}$ closely matches once more $g^\mathcal{B}$, with discrepancies concentrated near the discontinuity region, where the numerical mapping is least regular.
We omit here the results from the sampling analysis for the reflector obtained in the periodic setting as it was indistinguishable from \cref{subfig:heatmap_diff_exdisc}.

\section{Conclusion}

In this work, we introduced a mixed elliptic--hyperbolic Monge--Ampère formulation for a transport problem motivated by freeform optical design, and in particular, with the goal of finding differential periodic reflector surfaces. The proposed model couples elliptic and hyperbolic Monge--Ampère equations through a decomposition of the source domain, allowing the construction of transport mappings generated by potentials exhibiting a mixed curvature signature behavior. Through a discussion on boundary orientation (see Remark \ref{rem:boundary_orientation}) and boundary segmentation (see Proposition \ref{prop:boundary_segmentation_mix}), we derived Dirichlet boundary conditions which are useful in enforcing the mapping of the boundary of the target domain $\partial\calY$ to the boundary of the two source subdomains $\partial\calX^\tL$ and $\partial\calX^\tR$.
A least-squares method was developed for the numerical approximation of solutions to the resulting system. Numerical examples demonstrated the capability of the method to compute both isolated and periodic mixed-type solutions, inducing smooth reflector configurations relevant to optical array design. The mixed formulation introduced opens the possibility of studying a broader class of Monge--Ampère transport problems.
    
However, several questions remain open. From an analytical prespective, the existence and regularity properties of solutions to the mixed system are currently unknown. The current formulation is restricted to rectangular source and target domains, and extending the proposed framework to more general geometries constitutes an important next step. Moreover, the decomposition considered here partitions the source domain into two regions corresponding to elliptic and hyperbolic behavior. More general decompositions involving multiple subdomains, for example an $n\times m$ partition of the source domain, may allow for the construction of more complex reflector configurations.

\bibliographystyle{siam}
\bibliography{references}
\end{document}